\magnification=1100
\overfullrule0pt

\input amssym.def
\input prepictex
\input pictex
\input postpictex

%\voffset= .5 in

% ********************* Definitions ************************************

%\def\widetilde{\mathaccent"0365 }

\def\CC{{\Bbb C}}

\def\RR{{\Bbb R}}
\def\ZZ{{\Bbb Z}}

\def\cC{{\cal C}}

\def\Card{\hbox{Card}}

\def\codim{\hbox{codim}}

\def\diag{\hbox{diag}}
\def\id{\hbox{id}}
\def\im{\hbox{im}}

% ********************* FONTS ************************************

\font\smallcaps=cmcsc10
\font\titlefont=cmr10 scaled \magstep1

\font\sectionfont=cmbx10
\font\tinyrm=cmr10 at 8pt

% ******************** SECTION HEADERS ***************************

\newcount\sectno
\newcount\subsectno
\newcount\resultno

\def\section #1. #2\par{
\sectno=#1
\resultno=0
\bigskip\noindent{\sectionfont #1.  #2}~\medbreak}

\def\subsection #1\par{\bigskip\noindent{\it  #1} \medbreak}

%******************* MATHEMATICAL LABELS **************************

\def\prop{ \global\advance\resultno by 1
\medskip\noindent{\bf Proposition \the\sectno.\the\resultno. }\sl}
\def\lemma{ \global\advance\resultno by 1
\medskip\noindent{\bf Lemma \the\sectno.\the\resultno. }
\sl}
\def\fact{ \global\advance\resultno by 1
\medskip\noindent{\bf Fact \the\sectno.\the\resultno. }
\sl}

\def\remark{ \global\advance\resultno by 1
\medskip\noindent{\bf Remark \the\sectno.\the\resultno. }}
\def\example{ \global\advance\resultno by 1
\medskip\noindent{\bf Example \the\sectno.\the\resultno. }\sl}
\def\cor{ \global\advance\resultno by 1
\medskip\noindent{\bf Corollary \the\sectno.\the\resultno. }\sl}
\def\thm{ \global\advance\resultno by 1
\medskip\noindent{\bf Theorem \the\sectno.\the\resultno. }\sl}
\def\defn{ \global\advance\resultno by 1
\medskip\noindent{\it Definition \the\sectno.\the\resultno. }\slrm}

\def\endlemma{\rm\medskip}

\def\pf{\rm\smallskip\noindent{\it Proof. }}
\def\endpf{\qed\hfil\medskip}

%Homemade Struts:
\newbox\strutAbox
\setbox\strutAbox=\hbox{\vrule height 12pt depth6pt width0pt}
\def\strutA{\relax\copy\strutAbox}
\newbox\strutBbox
\setbox\strutBbox=\hbox{\vrule height 10pt depth5pt width0pt}
\def\strutB{\relax\copy\strutBbox}
\newbox\strutDbox
\setbox\strutDbox=\hbox{\vrule height 11pt depth5pt width0pt}

%high strut:
\newbox\strutHbox
\setbox\strutHbox=\hbox{\vrule height 11pt depth1pt width0pt}
\def\strutH{\relax\copy\strutHbox}
%low strut:
\newbox\strutLbox
\setbox\strutLbox=\hbox{\vrule height 1pt depth5pt width0pt}
\def\strutL{\relax\copy\strutLbox}

% hack to ignore lots of typed stuff....
\def\ignore#1{\relax}

% ******************  QED SIGNS  *********************************

\def\qed{\hbox{\hskip 1pt\vrule width4pt height 6pt depth1.5pt \hskip 1pt}}

\def\sqr#1#2{{\vcenter{\vbox{\hrule height.#2pt
\hbox{\vrule width.#2pt height#1pt \kern#1pt
\vrule width.2pt}
\hrule height.2pt}}}}

  % open square

%*************** EQUATIONS WITH NUMBERS **************

\def\formula{\global\advance\resultno by 1
\eqno{(\the\sectno.\the\resultno)}}
\def\formulano{\global\advance\resultno by 1 (\the\sectno.\the\resultno)}
\def\tableno{\global\advance\resultno by 1
\the\sectno.\the\resultno. }
\def\lformula{\global\advance\resultno by 1
\leqno(\the\sectno.\the\resultno)}

%********** DATING ******************************************
\def\monthname {\ifcase\month\or January\or February\or March\or April\or
May\or June\or
July\or August\or September\or October\or November\or December\fi}

\newcount\mins  \newcount\hours  \hours=\time \mins=\time
\def\now{\divide\hours by60 \multiply\hours by60 \advance\mins by-\hours
     \divide\hours by60         % NOTE: \divide only gives integer answers.
     \ifnum\hours>12 \advance\hours by-12
       \number\hours:\ifnum\mins<10 0\fi\number\mins\ P.M.\else
       \number\hours:\ifnum\mins<10 0\fi\number\mins\ A.M.\fi}

%**************** PAGE HEADERS *************************

\nopagenumbers
\def\runningtitle{\smallcaps graded hecke algebras}
\headline={\ifnum\pageno>1\eoheadline\else\firstheadline\fi}
\def\names{\smallcaps a. ram \enspace and \enspace a. shepler}
\def\firstheadline{}
\def\eoheadline{\ifodd\pageno\oddheadline\else\evenheadline\fi}
\def\oddheadline{\tenrm\hfil\runningtitle\hfil\folio}
\def\evenheadline{\tenrm\folio\hfil{\names}\hfil}

%**************** TITLE *************************
\vphantom{$ $}  %My kludge to get the first page to move down a bit
\vskip.75truein
\centerline{\titlefont Classification of graded Hecke algebras}
\smallbreak
\centerline{\titlefont for complex reflection groups}
\bigskip
\centerline{\rm Arun Ram}
%${}^\ast$ 
\centerline{Department of Mathematics}
\centerline{University of Wisconsin, Madison}
\centerline{Madison, WI 53706 USA}
\centerline{{\tt ram@math.wisc.edu}}
\medskip
\centerline{\rm and}
\medskip
\centerline{\rm Anne V. Shepler}
\centerline{Department of Mathematics}
\centerline{University of North Texas}
\centerline{Denton, TX 76203 USA}
\centerline{{\tt ashepler@unt.edu}}
%\medskip
%\centerline{Version: \today}

\footnote{}{\tinyrm 
Research of the first author supported in part by the National
Security Agency and by EPSRC Grant GR K99015 at the Newton Institute
for Mathematical Sciences.
Research of the second author supported in part by 
National Science Foundation grant DMS-9971099.}
\footnote{}{\tinyrm
}

\bigskip

%**************** ABSTRACT *************************
\noindent{\bf Abstract.}

The graded Hecke algebra for a finite Weyl group is intimately
related to the geometry of the Springer correspondence.
A construction of Drinfeld produces an analogue of a 
graded Hecke algebra for any finite subgroup of $GL(V)$.
This paper classifies all the algebras obtained by 
applying Drinfeld's construction to complex reflection 
groups.  By giving explicit (though nontrivial) isomorphisms,
we show that the graded Hecke algebras for 
finite real reflection groups constructed by 
Lusztig are all isomorphic to algebras obtained by Drinfeld's 
construction.  The classification shows that there
exist algebras obtained from Drinfeld's construction which are
not graded Hecke algebras as defined by Lusztig
for real as well as complex reflection groups.

%**************** KEY WORDS *************************
%\bigskip
%\noindent{\bf Keywords:} Iwahori-Hecke algebras, representations,
%symmetric functions
%%%%%%%%%%%%%%%%%%%%%%%%%%%%%%%%%%%%%%%%%%%%%%%%%%%%%%%%%%%

\section 0. Introduction

This paper is motivated by a general effort to generalize
the theory of Weyl groups and their
relation to groups of Lie type to the 
setting of complex reflection groups.   One natural
question is whether there are affine Hecke algebras
corresponding to complex reflection groups.  If they
exist then it might be possible to use these algebras
to build an analogue of the Springer correspondence 
for complex reflection groups.

A priori, one knows how to construct affine Hecke algebras
corresponding only to Weyl groups since both a finite real
reflection group $W$ and a $W$-invariant lattice (the
existence of which forces $W$ to be a Weyl group) 
are needed in the construction.  
Our search for analogues of graded Hecke algebras for 
complex reflection groups was motivated
by Lusztig's
results [Lu2] showing that the geometric information contained 
in the affine Hecke
algebra can be recovered from the corresponding graded
Hecke algebra.
Lusztig [Lu] defines the graded Hecke algebra for a finite
Weyl group $W$ with reflection representation $V$.  Let
$t_g$, $g\in W$, be a basis for the group algebra $\CC W$ of
$W$ and let $k_\alpha\in \CC$ be ``parameters'' indexed by the roots
in the root system of $W$ such that $k_\alpha$ depends
only on the length of the root $\alpha$. 
Then the graded Hecke algebra $H_{\rm gr}$ 
depending on the parameters $k_\alpha$ 
is the (unique) algebra structure on $S(V)\otimes \CC W$ 
such that
\smallskip\noindent
\itemitem{(a)}  the symmetric algebra of $V$,  $S(V)=S(V)\otimes 1$,
is a subalgebra of $H_{\rm gr}$,
\smallskip\noindent
\itemitem{(b)}  
the group algebra $\CC W= 1\otimes \CC W=\hbox{span-}\{1\otimes t_g\ |\ g\in W\}$ 
is a subalgebra of $H_{\rm gr}$, and
\smallskip\noindent
\itemitem{(c)} $t_{s_i} v = (s_i v)t_{s_i} - k_{\alpha_i}
\langle v,\alpha_i^\vee\rangle$
for all $v\in V$ and simple reflections $s_i$ in the simple roots
$\alpha_i$.  
\smallskip\noindent
This definition applies to all finite real reflection
groups $W$ since the simple roots and simple reflections are well
defined.
Unfortunately, the need for simple reflections
in the construction makes it unclear how to define
analogues for complex reflection groups.

For finite real reflection groups, the graded
Hecke algebra $H_{\rm gr}$ is a ``semidirect 
product'' of the polynomial ring $S(V)$ and the group
algebra $\CC W$.  Drinfeld [Dr] defines a {\it different}
type of semidirect product of $S(V)$ and $\CC W$, and 
Drinfeld's construction applies to all finite subgroups $G$
of $GL(V)$.  In this paper, we
\item{(1)} Classify all the algebras obtained by 
applying Drinfeld's construction to finite complex reflection
groups $G$, 
\item{(2)} Show that every graded Hecke algebra $H_{\rm gr}$
(as defined by Lusztig) for a finite real reflection group
is isomorphic to an algebra obtained by Drinfeld's
construction by giving explicit isomorphisms between these 
algebras.
\smallskip\noindent
The results from (2) show how Drinfeld's construction is a true
generalization of Lusztig's construction of graded Hecke algebras,
something which is not obvious.  Our classification
in (1) gives a complete solution to the problem of finding
all graded Hecke algebras for finite reflection groups.

A consequence of our classification is that there exist
graded Hecke algebras for finite real reflection groups which 
are not obtained with Lusztig's construction.
In this sense, Drinfeld's construction is a strict
generalization of the algebras $H_{\rm gr}$.  These
new algebras, and the algebras corresponding
to complex reflection groups that are not real reflection
groups, deserve further study and probably have interesting
representation theories.

For us, one surprising result of our classification 
is that no nontrivial graded Hecke algebra structures
exist for many complex reflection groups.
In some sense, this is disappointing, as we would have liked to find
nontrivial and interesting structures for {\it each} complex 
reflection group.

It might be that we have not yet hit upon the
``right'' definition of graded Hecke algebras.   For
example, we show that
there do not exist nontrivial graded Hecke algebra
structures, according to Drinfeld's definition,
for any of the complex reflection groups
$G(r,1,n)=(\ZZ/r\ZZ)\wr S_n$ when $r>2$ and $n>3$.
On the other hand, in the last section of this paper
we are able to construct algebras that ``look'' like they ought
to be graded Hecke algebras corresponding to these groups.
Is it possible that there is a ``better'' definition
of graded Hecke algebras which applies to complex reflection
groups and which includes the algebras that we introduce
in Section 5 as examples?

\medskip\noindent
{\bf Acknowledgements.}  We thank C. Kriloff for numerous stimulating
conversations during our work on this paper.  A. Ram thanks the Newton
Institute for the Mathematical Sciences at Cambridge University
for hospitality and support (EPSRC Grant No. GR K99015)
during Spring 2001 when the writing of this paper was completed.

\section 1. \ Graded Hecke algebras

In this section, we define the graded Hecke algebra following
Drinfeld [Dr].  Our main result in this section is Theorem 1.9c,
which determines exactly how many degrees of freedom one has
in defining a graded Hecke algebra.

Let $V$ be an $n$ dimensional vector space over $\CC$ 
and let $G$ be a finite subgroup of $GL(V)$.
The group algebra of $G$ is 
$$\CC G= \hbox{$\CC$-span}\{t_g \ |\ g \in G\},
\qquad\hbox{with}\quad t_g t_h = t_{gh}.$$
Let $a_g\colon V\times V \longrightarrow \CC$
be skew symmetric bilinear forms indexed by the elements 
of $G$ and let $A$ be the associative algebra generated 
by $V$ and $\CC G$ with the additional relations
$$
 t_hvt_{h^{-1}} = hv
 \qquad\hbox{and}\qquad
 [v,w] = \sum_{g\in G} a_g(v,w)\,  t_g,
 \qquad\hbox{
 for $h\in G$ and $v,w\in V$,}
\formula
$$
where $[v,w]=vw-wv$.
These relations allow every element $a\in A$ to be written
in the form
$$a = \sum_{g \in G} p_gt_g, \qquad p_g\in S(V),
\formula$$
where $S(V)$ is the symmetric algebra of $V$.  
More precisely, one must fix a section of the canonical 
surjection $T(V)\to S(V)$ from the tensor algebra of $V$ to $S(V)$
and take the elements $p_g$ to be in the image of this section.

The structure of $A$ depends on the choices of the
``parameters'' $a_g(v,w)\in \CC$.  
Our goal is to determine when the algebra $A$
will be a ``semidirect product'' of $S(V)$ and $\CC G$.  This
idea motivates the following definition [Dr].

\smallskip
The algebra $A$ is a {\it graded Hecke algebra} for $G$ if
$$A\cong S(V)\otimes \CC G
\quad \hbox{as a vector space,}$$
or, equivalently, if the expression in (1.2) is unique
for each $a\in A$.
A general element $a\in A$ is a linear combination of products of
elements $t_g$ and $u_i$, where $\{u_1,u_2,\ldots,u_n\}$ is a basis of $V$.  
There are 
two straightening operations needed to put $a$ in the form (1.2):
$$\hbox{(a)\ \ moving $t_h$'s to the right,
\qquad and \qquad
(b)\ \ putting $u_i u_j$ pairs in the correct order.}$$
These two straightening operations correspond to the two identities in (1.1).
Note that the ``correct order'' of $u_i u_j$ is determined by the choice 
of the section of the projection $T(V)\to S(V)$.  
%Let us assume that in the preferred ordering,
%$u_3$ is to the left of $u_2$ and $u_2$ is to the left of $u_1$. 
Let $v_1, v_2, v_3$ be arbitrary elements of $V$ and
let $h\in G$.
Applying the straightening operations to $t_hv_1v_2$ gives
$$\eqalign{
t_hv_1v_2 &= t_h[v_1,v_2]+t_hv_2v_1 \cr
&=t_h[v_1,v_2]+(hv_2)(hv_1)t_h \cr
}\qquad\qquad
\matrix{
\hbox{(rearrange $v_1$ and $v_2$)} \cr
\hbox{(move $t_h$ to the right),}\hfill \cr
}$$
and applying the straightening operations in a different order gives
$$\eqalign{
t_hv_1v_2 &= (hv_1)(hv_2)t_h \cr
&=[hv_1,hv_2]t_h+(hv_2)(hv_1)t_h. \cr}
\qquad\qquad 
\matrix{
\hbox{(move $t_h$ to the right)} \cr
\cr
}
$$
Setting these two expressions equal gives the relation
$$t_h[v_1,v_2]t_{h^{-1}}=[hv_1,hv_2],
\qquad\hbox{for all $h\in G$, $v_1,v_2\in V$}.\formula$$
Similarly, applying the straightening operations to
$v_1v_2v_3$ gives
$$\eqalign{
v_1v_2v_3
&=[v_1,v_2]v_3+v_2v_1v_3   \cr
&=[v_1,v_2]v_3+v_2[v_1,v_3]+v_2v_3v_1 \cr
&=[v_1,v_2]v_3+v_2[v_1,v_3]+[v_2,v_3]v_1+v_3v_2v_1  \cr
}
\qquad
\matrix{
\hbox{(moving $v_1$ past $v_2$)} \cr
\hbox{(moving $v_1$ past $v_3$)} \cr
\hbox{(straightening $v_2$ and $v_3$),} \cr
}$$
and applying the straightening operations in a different order gives
$$\eqalign{
v_1v_2v_3
&=v_1[v_2,v_3]+ v_1v_3v_2  \cr
&=v_1[v_2,v_3]+[v_1,v_3]v_2+v_3v_1v_2 \cr
&=v_1[v_2,v_3]+[v_1,v_3]v_2+v_3[v_1,v_2]+v_3v_2v_1 \cr
}\qquad
\matrix{
\hbox{(moving $v_2$ past $v_3$)} \cr
\hbox{(moving $v_1$ past $v_3$)} \cr
\hbox{(straightening $v_1$ and $v_2$).} \cr
}
$$
These are equal if
$$[v_1,[v_2,v_3]]+[v_2,[v_3,v_1]]+[v_3,[v_1,v_2]]=0,
\qquad\hbox{for all $v_1,v_2,v_3\in V$.}
\formula$$
Conversely, the identities (1.3) and (1.4) are exactly what is
needed to guarantee that any order of application of the
straightening operations (a) and (b) will produce the 
same normal form (1.2) for the element $a$.  Thus
we have

\lemma  Let $A$ be an algebra defined as in (1.1).  Then
$A$ is a graded Hecke algebra if and only if the identities
(1.3) and (1.4) hold in $A$.
\endlemma

Using (1.1), the relations (1.3) and (1.4) 
can be rewritten in terms of the bilinear forms 
$a_g\colon V\times V\to \CC$ as
$$
a_g( v_1, v_2) = a_{hgh^{-1}}(hv_1,hv_2) \qquad \hbox{ and }
\formula
$$
$$
a_g(v_3,v_1)(gv_2-v_2)+a_g(v_2,v_3)(gv_1-v_1)+a_g(v_1,v_2)(gv_3-v_3)=0
\formula
$$
for $v_1,v_2,v_3\in V$ and $g,h\in G$.

Let $\langle \ \ ,\  \,\rangle : V \times V \rightarrow \CC$
be a $G$-invariant nondegenerate Hermitian form on $V$.  
For each $g\in G$, set 
$$\eqalign{
V^g &= \{ v\in V\ |\ gv=v\}, \cr
(V^g)^\perp &= \{ v\in V\ |\ \hbox{$\langle v,w\rangle = 0$ for all $w\in V^g$}\}, 
\qquad\hbox{and}\cr
\ker a_g &= \{v\in V\ |\ \hbox{$a_g(v,w)=0$ for all $w\in V$}\}. \cr}
$$
\lemma  
Let $G$ be a finite subgroup of $GL(V)$ and let $g\in G$.
\smallskip\noindent
\item{(a)} $(V^g)^\perp = \{ v-gv\ |\ v\in V\}$.
\smallskip\noindent
\item{(b)} Suppose $g\ne 1$.  If $\codim(V^g)=2$ and $a\colon V\times V\to \CC$ is a skew
symmetric bilinear form such that $\ker a = V^g$, then $a$ satisfies (1.7).
\smallskip\noindent
Let $A$ be a graded Hecke algebra for $G$ defined by skew symmetric bilinear
forms $a_g\colon V\times V\to \CC$. 
\item{(c)} If $g\ne 1$ then $\ker a_g\supseteq V^g$.
\smallskip\noindent
\item{(d)} If $g\ne 1$ and $a_g\ne 0$ then $\ker a_g=V^g$ and $\codim(V^g)=2$.
\smallskip\noindent
\item{(e)} If $g\ne 1$ and $a_g\ne 0$ then, for all $h\in G$,
$$a_{h^{-1}gh}(b_1,b_2) = \det(h^\perp)a_g(b_1,b_2),$$
where $\{b_1,b_2\}$ is a basis of $(V^g)^\perp$ and $h^\perp\colon (V^g)^\perp
\to (V^g)^\perp$ is the composition of $h$ restricted to $(V^g)^\perp$
with the canonical projection $V\to V/V^g$.
\pf
(a) Consider the map $\phi\colon V\to V$ given by $\phi(v)=v-gv$.
Then $\ker \phi=V^g$ and $\im\, \phi\subseteq (V^g)^\perp$
since, if $v\in V, w\in V^g$, then
$$\langle v-gv,w\rangle=\langle v,w\rangle-\langle gv,w\rangle
=\langle v,w\rangle-\langle gv,gw\rangle
=\langle v,w\rangle-\langle v,w\rangle.$$
Since $\dim(\im\,\phi)={\rm codim}(\ker\phi)={\rm codim}(V^g)$
it follows that $\im\,\phi = (V^g)^\perp$.
%So $\langle v-gv,w\rangle=0$ for all $v \in V$
%if and only if
%$\langle v,w - g^{-1} w\rangle =0$
%for all $v \in V$ if and only if
%$w \in V^g$.
\smallskip\noindent
(b) Let $v_1, v_2, v_3 \in V$.  If any 
$v_i \in V^g$, then (1.7) holds trivially for the skew symmetric
form $a$.  So assume each $v_i \not\in V^g$ and write each $v_i$ as 
$v_i^+ + v_i^-$ where $v_i^+ \in V^g$ and $v_i^- \in (V^g)^\perp$.
Then
$$
 a(v_i,v_j) = a(v_i^-,v_j^-)
 \qquad\hbox{and}\qquad
 v_i -gv_i = v_i^- - gv_i^-.
$$
Since $\dim(V^g)^\perp = 2$, at least one of the $v_i^-$ is a 
linear combination of the other two.
Say $v_1^- = c_2v_2^- + c_3 v_3^-$ with $c_2$, $c_3 \in \CC$.
Substituting  $v_i -gv_i = v_i^- - gv_i^-$ and 
$v_1^- = c_2v_2^- + c_3 v_3^-$ then yields
$$
 \eqalign{
 a(v_3,v_1)&(gv_2-v_2)+a(v_2,v_3)(gv_1-v_1)+a(v_1,v_2)(gv_3-v_3) \cr
 &=a(v_3^-,v_1^-)(gv_2^--v_2^-)+a(v_2^-,v_3^-)(gv_1^--v_1^-)
 +a(v_1^-,v_2^-)(gv_3^--v_3^-) 
%\cr
% &=c_2a(v_3^-,v_2^-)(gv_2^--v_2^-)
% +c_3a(v_3^-,v_3^-)(gv_2^--v_2^-) \cr
% &\qquad+c_2a(v_2^-,v_3^-)(gv_2^--v_2^-)
% +c_3a(v_2^-,v_3^-)(gv_3^--v_3^-) \cr
% &\qquad+c_2a(v_2^-,v_2^-)(gv_3^--v_3^-) 
% +c_3a(v_3^-,v_2^-)(gv_3^--v_3^-) 
 =0, \cr
 }
$$
and so (1.7) holds.
\smallskip\noindent
(c)
Let $v_3 \in V^g$ and $v_2\in V$.
\smallskip\noindent
If $v_2 \in V^g$, then
$a_g(v_2,v_3)(gv_1-v_1)=0$ for all $v_1 \in V$
by (1.7).  Since $V^g\ne V$, there is some $v_1$ 
such that $gv_1\ne v_1$ and so $a_g(v_2,v_3)=0$.
\smallskip\noindent
If $v_2 \not\in V^g$, let 
$\displaystyle{v_1 = \sum_{k=1}^{r-1} {g^k}v_2}$,
where $r$ is the order of $g$.
By (1.6),
$a_g(v_3, g^{k} v_2) = a_g(g^{-k} v_3, v_2) = a_g(v_3,v_2)$, for any $k$,
and so
$$\eqalign{
0 &= a_g(v_3,v_1)(g v_2 -v_2) + a_g(v_2,v_3)(gv_1-v_1)  \cr
&=(r-1)a_g(v_3,v_2)(gv_2-v_2) + a_g(v_3,v_2)(gv_2-v_2) 
=r a_g(v_3,v_2)(gv_2 -v_2).\cr
}
$$
Thus $a_g(v_3,v_2) = 0$.
\smallskip\noindent
Hence, for all $v_2 \in V$, $a_g(v_3,v_2) = 0$ and so $V^g \subseteq \ker a_g$.
\smallskip\noindent
(d)
By (c), $\codim(V^g)\ge \codim(\ker a_g)$.
Since $a_g \neq 0$, there exist $v, w \in V$ with
$a_g(v,w) \neq 0$ and so $\codim(\ker a_g) \geq 2$.
Let $v_1-gv_1$ and $v_2-gv_2$ be linearly independent elements of $(V^g)^\perp$.
Then (1.7) implies that any element $v_3 - gv_3 \in (V^g)^\perp$
is a linear combination of $v_1 - gv_1$ and $v_2 - gv_2$, and so 
$$2 \ge \dim((V^g)^\perp) = \codim(V^g) \ge \codim(\ker a_g) \ge 2.$$
Thus $V^g = \ker a_g$ and $\codim(V^g)=2$.
\smallskip\noindent
(e)  Write 
$hb_1 = h_{11}b_1+h_{21}b_2+(hb_1)^g$ and
$hb_2 = h_{12}b_1+h_{22}b_2+(hb_2)^g$ 
with $h_{ij}\in \CC$ and $(hb_i)^g\in V^g$.  Then
$$\eqalign{
a_{h^{-1}gh}(b_1,b_2)&=a_g(hb_1,hb_2)
=a_g(h_{11}b_1+h_{21}b_2+(hb_1)^g,
h_{12}b_1+h_{22}b_2+(hb_2)^g) \cr
&=(h_{11}h_{22}-h_{21}h_{12})a_g(b_1,b_2)
=\det(h^\perp)a_g(b_1,b_2) \cr}$$
since $a_g$ is skew symmetric and $V^g\subseteq \ker a_g$.
\endpf

The following theorem is a slightly strengthened
version of statements (given without proof) in [Dr].

\thm Let $G$ be a finite subgroup of $GL(V)$ and let
$Z_G(g)=\{ h\in G\ |\ hg=gh\}$ denote the centralizer of an element $g$ in $G$.
\item{(a)} If $A$ is a graded Hecke algebra for $G$,
then the values of $a_{h^{-1}gh}$ are determined by 
the values of $a_g$ via the equation
$$a_{h^{-1}gh}(v_1,v_2) = a_g(hv_1,hv_2),
\qquad\hbox{for all $g,h\in G$, $v_1,v_2\in V$}.$$
\item{(b)}  
For $g\ne 1$, there is a graded
Hecke algebra $A$ with $a_g\ne 0$ if and only if
$$\ker a_g=V^g,
\quad \codim(V^g)=2,   \quad\hbox{and}\quad   \det(h^\perp)=1,
\enspace\hbox{for all $h\in Z_G(g)$,}
$$
where $h^\perp$ is $h$ restricted to the space $(V^g)^\perp$.
In this case, $a_g$ is determined by its value $a_g(b_1,b_2)$ on a basis
$\{b_1,b_2\}$ of $(V^g)^\perp$.
\item{(c)} 
Let $d$ be the number of conjugacy classes of $g\in G$ such that 
$\codim(V^g)=2$ and $\det(h^\perp)=1$ for all $h\in Z_G(g)$,
where $h^\perp$ is $h$ restricted to the space $(V^g)^\perp$.
The sets $\{a_g\}_{g\in G}$ corresponding
to graded Hecke algebras $A$ form a vector space of dimension 
$d+\dim((\bigwedge^2 V)^G)$.
\pf (a) is simply a restatement of (1.6).  
\smallskip\noindent
(b) $\Longrightarrow$:  If $A$ is a graded Hecke algebra and $a_g\ne 0$ then
by Lemma 1.8d, $\codim(V^g) = 2$ and $\ker a_g = V^g$.
So $a_g$ is determined by its value $a_g(b_1,b_2)$ 
on a basis $b_1,b_2$ of $(V^g)^\perp$.  Suppose $h \in Z_G(g)$.
Then, by Lemma 1.8e,
$$ a_g(b_1, b_2) =
 a_{hgh^{-1}}(hb_1, hb_2) = a_g(h b_1, h b_2) =  \det(h^\perp) \ a_g(b_1, b_2),
$$
and so $\det(h^\perp) = 1$.
Note that 
$h(V^g) = V^g$ and $h(V^g)^\perp=(V^g)^\perp$ 
since, for  each $v \in V^g$, $h(v) = hg(v)= gh(v)$.  

\smallskip\noindent
$\Longleftarrow$:  
If $\codim(V^g)=2$ then, up to constant multiples,
there is a unique skew symmetric form on $V$ which is
nondegenerate on $(V^g)^\perp$ and which has $\ker a_g=V^g$.  
Fix such a form and then define forms $a_h$, $h\in G$, by
$$a_h(v_1,v_2) = \cases{
a_g(kv_1,kv_2) & if $h=k^{-1}gk$, \cr
0 &otherwise, \cr}
\formula$$
for $v_1,v_2\in V$.  Let 
$a_1$ be any $G$-invariant skew symmetric form on $V$.
Then this collection $\{a_g\}_{g\in G}$ of skew symmetric bilinear
forms satisfies (1.6) by definition
and (1.7) by Lemma 1.8b.  Thus (by Lemma 1.5), 
it determines a graded Hecke algebra $A$ via (1.1).  

\smallskip\noindent
(c) From (a) and (b) it follows that
the sets $\{a_g\}_{g\in G}$, running over all graded Hecke algebras $A$ 
for $G$, form a vector space.  
Since each of the collections $\{a_g\}_{g\ne 1}$ constructed by (1.10)
has its support on a single conjugacy class, these collections
form a basis of the vector space of sets $\{a_g\}_{g\ne 1}$.
The only condition on the form $a_1$ is that it
satisfies (1.6), which means that it is a $G$-invariant
element of $(\bigwedge^2 V)^*$.
\endpf

The following consequence of Theorem 1.9 will be useful for completing the
classification of graded Hecke algebras for complex reflection
groups.

\cor  Assume that $G$ contains $h=\xi\cdot 1$ for
some $\xi \in \CC \, \backslash \, \{\pm1\}$.   If $A$ is a graded Hecke algebra for $G$,
then $a_g=0$ for all $g\ne 1$. 
\pf
If $h = \xi\cdot 1 \in G$, then $h\in Z_G(g)$ for every $g\in G$ and 
$\det(h^\perp) = \xi^2$ if $\codim(V^g)=2$.  The statement then
follows from Theorem 1.9b.
\endpf

\vfill\eject

\section 2. The classification for reflection groups 

A {\it reflection} is an element of $GL(V)$ that has exactly
one eigenvalue not equal to $1$.  The {\it reflecting hyperplane}
of a reflection is the $(n-1)$-dimensional subspace which is fixed
pointwise. A {\it complex reflection group} $G$
is a finite subgroup of $GL(V)$ generated by reflections.
The group $G$ is {\it irreducible} if $V$ cannot be written in the form
$V = V_1\oplus V_2$ where $V_1$ and $V_2$ are $G$-invariant
subspaces.  The group $G$ is a {\it real reflection group} if 
$V=\CC\otimes_\RR V_\RR$ for a real vector space $V_\RR$ 
and $G\subseteq GL(V_\RR)$.

The following facts about reflection groups are well known.

\lemma Let $G$ be an irreducible reflection group.
\itemitem{(a)}  [ST, Theorem 5.3] The number of elements $g\in G$ such that
$\codim(V^g)=2$ is $\sum_{i<j} m_i m_j$ where $m_1,\ldots, m_n$
are the exponents of $G$.
\itemitem{(b)} [Ca, Lemma 2]  If $G$ is a real reflection group and 
$g\in G$ with $\codim(V^g)=2$, then $g$ is
the product of two reflections.
\itemitem{(c)} [OT, Theorem 6.27] For any $g\in G$, the space $V^g$ is the 
intersection of reflecting hyperplanes.
\endlemma

\noindent
{\bf Remark.}  The statement of Lemma 2.1b does not hold for complex
reflection groups.  Consider the exceptional complex reflection
group $G_4$ of rank 2, in the notation of Shephard and Todd [ST].
All the reflections have order $3$ and $-1\in G_4$.  Suppose
$-1=rs$ for two reflections $r$ and $s$.  
If $s$ has eigenvalues $1$ and $\omega$, where $\omega$ is a primitive
cube root of unity, then $r^{-1}=-s$ has eigenvalues
$-1$ and $-\omega$, a contradiction to the assumption that $r$ is a 
reflection.  Thus $-1\in G_4$ is not a product of two reflections.

\lemma
Let $G\subseteq GL(V)$ be a complex reflection group.
 Let $A$ be a
graded Hecke algebra for $G$ and let $g\in G$.
Let $V^G=\{v\in V\ |\ \hbox{$gv=v$ for all $g\in G$}\}$
be the invariants in $V$. 
\item{(a)}  If $g=1$  and $\dim V^G \le 1$, then $a_g=0$.
\item{(b)}  If the order of $g$ is $2$, then $a_g=0$.
\pf
(a)  
Let $\langle \, ,\rangle\colon V\times V\to \CC$ be a nondegenerate $G$-invariant
Hermitian form on $V$ and write $V = V^G\oplus (V^G)^\perp$ where 
$(V^G)^\perp = \{v\in V\ |\ \hbox{$\langle v,w\rangle = 0$ for all $w\in V^G$}\}$.
Since $\dim(V^G)\le 1$ and $a_1$ is skew symmetric, $a_1$ restricted 
to $V^G$ is $0$.  There is a basis $\alpha_1,\ldots,\alpha_k$ of 
$(V^G)^\perp$ and constants $\xi_1,\ldots,\xi_k\in \CC$, $\xi_i \neq 1$,
such that the reflections $s_1,\ldots, s_k$ given by
$$
 s_iv = v+(\xi_i-1){\langle v,\alpha_i\rangle\over 
 \langle \alpha_i,\alpha_i\rangle} \alpha_i,
 \qquad\hbox{for $v\in V$},
$$
are in $G$.  Equation (1.6)
implies that, for any $v\in V$,  
$$
 a_1(\alpha_i, v) = a_1(s_i\alpha_i, s_iv)
 =a_1\left(\xi_i \alpha_i, v +(\xi_i-1) {\langle v,\alpha_i \rangle
       \over \langle \alpha_i, \alpha_i \rangle}   \alpha_i\right)
 =\xi_i a_1(\alpha_i,v),\
$$
since $a_1(\alpha_i, \alpha_i) = 0$
(as $a_1$ is skew symmetric).  Since $\xi_i\ne 1$, 
$a_1(\alpha_i, v)=0$ for $1\leq i \leq k$.
Thus $\ker a_1=V$.
\smallskip\noindent
{(b)}
Since $g^2=1$, all eigenvalues of $g$ are $\pm 1$.
If $\codim(V^g) \neq 2$,
then $a_g = 0$ by Theorem 1.9b.
If $\codim(V^g)=2$, then
$$
 g = \id_{V^g} \oplus (-\id_{(V^g)^\perp})
$$
as a linear transformation on $V$.  By [St1, Theorem 1.5],
[Bou V, \S 5 Ex.~8],
the stabilizer, ${\rm Stab}(V^g)$, of $V^g$
is a reflection subgroup of $G$ and so there is a reflection
$s\in {\rm Stab}(V^g)$ that is the identity on $V^g$.
So $s\in Z_G(g)$ and $\det(s) = \det(s^\perp)\ne 1$, where
$s^\perp$ is $s$ restricted to $(V^g)^\perp$. 
Thus, by Theorem 1.9b, $a_g=0$.
\endpf

%Using the classification of finite subgroups of the unitary
%group $\hbox{U}_2(\CC)$, one can show
%that $Z_G(g)$ is cyclic whenever there is a graded
%Hecke algebra for $G$ with $a_g \neq 0$.

\vfill\eject

\subsection 2A.\ \ Real reflection groups

If $G\subseteq GL(V)$ is a real reflection group 
then $V=\CC\otimes_\RR V_\RR$
and $G\subseteq GL(V_\RR)$, where $V_\RR$ is a real vector space.
We shall assume that $G$ is irreducible.

Let us recall some basic facts about real reflection groups
which can be found in [Hu] or [Bou].  The action of $G$ on $V_\RR$ has
fundamental chambers $wC$ indexed by $w\in G$.  The {\it roots}
for $G$ are vectors $\alpha\in V_\RR$ such that the reflections
in $G$ are the reflections $s_\alpha$ in the hyperplanes
$$H_\alpha = \{ v\in V_\RR\ |\ \langle v,\alpha\rangle=0\}.$$
For each fundamental chamber $C$, the reflections
$s_1,s_2,\ldots,s_n$ in the hyperplanes $H_{\alpha_1},
H_{\alpha_2},\ldots, H_{\alpha_n}$ that bound $C$ form a set of 
{\it simple reflections} for $G$.  The simple reflections obtained
from a different choice of fundamental chamber $wC$ are 
$ws_1w^{-1},\ldots, ws_nw^{-1}$.

\thm
Let $G\subseteq GL(V_\RR)$ be a real reflection group.
Let $s_1,\ldots, s_n$ be a set of simple reflections
in $G$ and let $m_{ij}$ be the order of $s_is_j$.
Then $g\in G$ satisfies $g^2\ne 1$,
$\codim(V^g)=2$, and $\det(h^\perp)=1$ for all $h\in Z_G(g)$ 
(the conditions in Theorem 1.9c) if and only if $g$ is conjugate to 
$$
 (s_is_j)^k,  \qquad \hbox{with $0< k< { m_{ij} \over  2}$,}
$$
for some $1\le i,j\le n$.
\pf
$\Longrightarrow:$  Let $\alpha$ and $\beta$ be two roots such that 
$V^g=H_{\alpha}\cap H_{\beta}$ (see Lemma 2.1c).
Then $H_{\alpha}\cap H_{\beta}$ has nontrivial
intersection with some fundamental chamber $C$ for $W$,
and we may assume that $H_{\alpha}$ and $H_{\beta}$ are
walls of the chamber $C$ (since $C$ is a cone in $\RR^n$).
Since choosing simple reflections with respect to a 
different chamber $wC$ corresponds to conjugation 
by $w$, we may assume that the reflections 
in the hyperplanes
$H_{\alpha}$ and $H_{\beta}$ are simple reflections
and $\alpha = \alpha_1$ and $\beta = \alpha_2$.

The element $g$ is an element of the stabilizer
${\rm Stab}(V^g)$, which is a reflection
group by [St1, Theorem 1.5].  Since $\codim(V^g)=2$, ${\rm Stab}(V^g)$ is
a rank two real reflection group, and therefore 
a dihedral group.  This dihedral group
is generated by the two simple reflections $s_1$ and $s_2$ in the 
hyperplanes $H_{\alpha_1}$ and $H_{\alpha_2}$
(restricted to $(V^g)^\perp$) and all reflections have determinant $-1$.   
Let $g^\perp$ be the element $g$ restricted to $(V^g)^\perp$. 
Since $g\in Z_G(g)$, $\det(g^\perp)=1$, and so $g$ must
be a product of an even number of reflections.
Thus $g=(s_1s_2)^k$ or $g=(s_2s_1)^k$, for some $0<k\le m/2$,
where $m$ is the order of $s_1s_2$.  Since $g^2\ne 1$, $k\ne m/2$, 
and so $g$ is conjugate to $(s_1s_2)^k$ for some $0<k<m/2$.

\smallskip\noindent
$\Longleftarrow:$  Assume that $g=(s_is_j)^k$ for some
$0<k<m_{ij}/2$.  Then $V^g = H_{\alpha_i}\cap H_{\alpha_j}$
and so $\codim(V^g)=2$.  
Since $g$ is a product of an even number of reflections, 
$\det(g^\perp)=1$.  The only elements of $O(V_\RR)\cong O_2(\RR)$
that are diagonalizable in $GL(V_\RR)\cong GL_2(\RR)$ are
$\pm 1$ and elements with determinant $-1$.
Thus, the eigenvectors of the element $g^\perp$
(which has distinct eigenvalues since it is not $\pm 1$)
do not lie in $V_\RR$, only in $V=\CC\otimes_\RR \RR$.
Let $h\in Z_G(g)$ and let $h^\perp \in O(V_\RR)\cong O_2(\RR)$ 
denote $h$ restricted to 
$(V^g)^\perp$.  Since $h^\perp$ commutes with $g^\perp$
and $g^\perp$ has distinct eigenvalues,
$g^\perp$ and $h^\perp$ have the same eigenvectors.
Hence, $\det h^\perp = 1$.
\endpf

Using Theorem 2.3 and Theorem 1.9b, we can read off the graded Hecke algebras
for the irreducible real reflection groups from
the Dynkin diagrams.  For each irreducible real reflection group,
label a set of simple reflections 
$s_1,\ldots,s_n$ using the Dynkin diagrams below.
If nodes $i$ and $j$ and nodes $j$ and $k$ are connected by single
edges, then $s_is_j$ is conjugate to $s_js_k$ via the element $s_is_js_k$.

The following table gives representatives of the conjugacy classes
of $g\in G$ that may have $a_g\ne 0$ for some graded Hecke algebra $A$.
We assume that the reflection group $G$
is acting on its irreducible reflection representation 
$V$.  When $G$ is the symmetric group $S_n$ acting
on an $n$-dimensional vector space $V$ by permutation matrices, 
then $\dim(V^G) = 1$ and, by Lemma 2.2a and Theorem 2.3,
$a_g\ne 0$ for some graded Hecke algebra $A$ only if
$g$ is conjugate to the three cycle $(1,2,3)=s_1s_2$
(this example is analyzed in Section 3).

$$
{
\vbox{
\offinterlineskip
\hrule
\halign{& \vrule# & 
                \quad $#$ \hfil\quad
        & \vrule# &  \quad\hfil $#$ \hfil\quad 
        & \vrule#  \cr
 & \strutH \ {\rm Type} && {\rm Representative }\  g & \cr
 & \strutL     && {\rm with }\ a_g\, \ne\,  0 & \cr
\noalign{\hrule}
&  \strutB  A_{n-1}  && s_1s_2 &\cr
&  \strutB B_n && s_1s_2, \ s_2s_3 &\cr
&  \strutB D_n && s_2s_3 &\cr
&  \strutB E_6, E_7, E_8 && s_1 s_4 &\cr
&  \strutB F_4 && s_1s_2, \ s_2s_3, \ s_3s_4 &\cr
&  \strutB H_3, H_4 && s_1s_2, \ (s_1s_2)^2, \ s_2s_3 &\cr
&  \strutB I_2(m) && (s_1s_2)^k, \ \  0<k<m/2 &\cr
}
\hrule}}
$$
\centerline{{\bf Table 1.} Graded Hecke algebras for real reflection groups.}

\font\scriptsize=cmr10 at 5pt  %for the pictures

\vskip .3in
$$
%********************************************************************
% Dynkin diagram of type A_n
%********************************************************************
\beginpicture
\font\scriptsize=cmr10 at 5pt  %for the pictures
\setcoordinatesystem units <1cm,1cm>                % sets scale
\setplotarea x from 1.5 to 8.5, y from 1.5 to 2.5   % sets plot size up
\put {$A_{n-1}$} at 2 2.5
{\scriptsize
\multiput {$\circ$} at 3   2 *1 1 0 /      %
\multiput {$\circ$} at 7   2 *1 1 0 /      %  puts nodes in
\put {$^1$}   at 3 2.2   %
\put {$^2$}   at 4 2.2   %
\put {$^{n-2}$} at 7 2.2   % 
\put {$^{n-1}$}   at 8 2.2          % label nodes with roots
\linethickness=0.75pt           % sets line thickness
\putrule from 3.05 2 to 3.95 2  %
\putrule from 7.05 2 to 7.95 2  % puts solid lines between nodes
\setdashes <2mm,1mm>            %
\putrule from 4.05 2 to 6.95 2  % draws dotted lines between nodes
}
%\put{Dynkin diagram of type $A_{n-1}$}[b] at 5.5 1.25
\endpicture
\qquad
%********************************************************************
% Dynkin diagram of type B_n
%********************************************************************
\beginpicture
\setcoordinatesystem units <1cm,1cm>       % sets scale
\setplotarea x from 1.5 to 8.5, y from 1.5 to 2.5 % sets plot size up
\put {$B_n$} at 2 2.5
{\scriptsize
\multiput {$\circ$} at 3   2 *2 1 0 /      %
\multiput {$\circ$} at 7   2 *1 1 0 /      %  puts nodes in
\put {$^1$}     at 3 2.2   %
\put {$^2$}     at 4 2.2   %
\put {$^3$}     at 5 2.2   %
\put {$^{n-1}$}   at 7 2.2   %
\put {$^n$}     at 8 2.2   % label nodes with roots above
\linethickness=0.75pt                      % sets line thickness
\putrule from 3.03 2.045 to 3.97 2.045  %
\putrule from 3.03 1.955 to 3.97 1.955  %
\putrule from 4.05 2 to 4.95 2         %   puts solid lines between nodes
\putrule from 7.05 2 to 7.95 2         %
\setdashes <2mm,1mm>            %
\putrule from 5.05 2 to 6.95 2  % draws dotted lines between nodes
}
%\put{Undirected Dynkin diagram of type $C_n$}[b] at 5.5 1.25
\endpicture
$$
$$
%********************************************************************
% Dynkin diagram of type D_n
%********************************************************************
\beginpicture
\setcoordinatesystem units <1cm,1cm>                % sets scale
\setplotarea x from 1.5 to 8.5, y from 1 to 3           % sets plot size
\put {$D_n$} at 2 2.5 
{\scriptsize
\multiput {$\circ$} at 8   2 *1 -1 0 /      %
\multiput {$\circ$} at 5   2 *1 -1 0 /      %  puts nodes in
\multiput {$\circ$} at 3 1.5 *1 0 1 /      %
\put {$^n$}   at 8 2.2   %
\put {$^{n-1}$} at 7 2.2   %
\put {$^4$}   at 5 2.2   %
\put {$^3$}   at 4 2.2   %
\put {$^2$}   at 3 2.7   %
\put {$^1$}   at 3 1.2   % label nodes with roots below
\setplotsymbol ({\rm .})
\linethickness=0.75pt                           % sets line thickness
\putrule from 7.05 2 to 7.95 2       % puts solid lines between nodes
\putrule from 4.05 2 to 4.95 2       %
\setlinear
\plot 3.05 2.5   3.95 2.05 / %
\plot 3.05 1.5   3.95 1.95 / %
\setdashes <2mm,1mm>          %
\putrule from 5.05 2 to 6.95 2  % draws dotted lines between nodes
}
%\put{Dynkin diagram of type $D_n$}[b] at 5.5 0.75
\endpicture
\qquad
\beginpicture
\setcoordinatesystem units <1cm,1cm>                % sets scale
%********************************************************************
% Dynkin diagram of E_6
%********************************************************************
\setplotarea x from 1.5 to 8.5, y from 1 to 3  % sets plot size
\put {$E_6$} at 2 2.5
{\scriptsize
\multiput {$\circ$} at 7   2 *4 -1 0 /      % puts nodes in
\put{$\circ$} at 5   1.1     %
%\put{$\circ$} at 5   0.2     %
\put {$^6$}   at 7 2.2   %
\put {$^5$}   at 6 2.2   %
\put {$^4$}   at 5 2.2   % label nodes
\put {$^3$}   at 4 2.2   %
\put {$^2$}   at 3 2.2   % 
\put {$^1$}   at 4.8 1.05   % 
%\put {$^0$}   at 4.8 0.15  % 
\setplotsymbol ({\rm .})
\linethickness=0.75pt                           % sets line thickness
\putrule from 6.05 2 to 6.95 2       % 
\putrule from 5.05 2 to 5.95 2       %
\putrule from 4.05 2 to 4.95 2       %  puts solid lines between nodes
\putrule from 3.05 2 to 3.95 2       %
\putrule from 5 1.15 to 5 1.95       %
%\putrule from 5 -1.25 to 5 -.45     %
}
% end of dynkin diagram of E_6
\endpicture
$$
$$
\beginpicture
\setcoordinatesystem units <1cm,1cm>                % sets scale
%********************************************************************
% Dynkin diagram of E_7
%********************************************************************
\setplotarea x from 1 to 7.2, y from -1 to 1.5  % sets plot size
\put{$E_7$} at 1 1
{\scriptsize
\multiput {$\circ$} at 7   0.5 *5 -1 0 /      % puts nodes in
\put{$\circ$} at 4 -0.4     %
%\put{$\circ$} at 5 -1.3     %
\put {$^7$}   at 7 0.7   %
\put {$^6$}   at 6 0.7   %
\put {$^5$}   at 5 0.7   %
\put {$^4$}   at 4 0.7   % label nodes
\put {$^3$}   at 3 0.7   %
\put {$^2$}   at 2 0.7   % 
\put {$^1$}   at 3.8 -0.45   % 
%\put {$^0$}   at 4.8 -1.35   % 
\setplotsymbol ({\rm .})
\linethickness=0.75pt                           % sets line thickness
\putrule from 6.05 0.5 to 6.95 0.5       % 
\putrule from 5.05 0.5 to 5.95 0.5       %
\putrule from 4.05 0.5 to 4.95 0.5       %  puts solid lines between nodes
\putrule from 3.05 0.5 to 3.95 0.5       %
\putrule from 2.05 0.5 to 2.95 0.5       %
\putrule from 4 -.35 to 4 .45       %
%\putrule from 5 -1.25 to 5 -.45     %
}
% end of dynkin diagram of E_7
\endpicture
\qquad
\beginpicture
\setcoordinatesystem units <1cm,1cm>                % sets scale
%********************************************************************
% Dynkin diagram of extended E_8
%********************************************************************
\setplotarea x from 0 to 7.2, y from -1 to 1.5  % sets plot size
\put{$E_8$} at 0 1.2
{\scriptsize
\multiput {$\circ$} at 7   0.5 *6 -1 0 /      % puts nodes in
\put{$\circ$} at 3 -0.4     %
%\put{$\circ$} at 5 -1.3     %
\put {$^8$}   at 7 0.7   %
\put {$^7$}   at 6 0.7   %
\put {$^6$}   at 5 0.7   %
\put {$^5$}   at 4 0.7   % label nodes
\put {$^4$}   at 3 0.7   %
\put {$^3$}   at 2 0.7   % 
\put {$^1$}   at 2.8 -0.45   % 
\put {$^2$}   at 1 0.7   % 
\setplotsymbol ({\rm .})
\linethickness=0.75pt                           % sets line thickness
\putrule from 6.05 0.5 to 6.95 0.5       % 
\putrule from 5.05 0.5 to 5.95 0.5       %
\putrule from 4.05 0.5 to 4.95 0.5       %  puts solid lines between nodes
\putrule from 3.05 0.5 to 3.95 0.5       %
\putrule from 2.05 0.5 to 2.95 0.5       %
\putrule from 1.05 0.5 to 1.95 0.5       %
\putrule from 3 -.35 to 3 .45       %
%\putrule from 5 -1.25 to 5 -.45     %
}
% end of dynkin diagram of E_8
\endpicture
$$
$$
%********************************************************************
% Dynkin diagram of type F_4
%********************************************************************
\beginpicture
\setcoordinatesystem units <1cm,1cm>       % sets scale
\setplotarea x from 2.5 to 7, y from 1.5 to 2.8 % sets plot size up
\put{$F_4$} at 2 2.5
{\scriptsize
\multiput {$\circ$} at 3   2 *3 1 0 /      %  puts nodes in
%\put {$^0$}     at 7 2.2   %
\put {$^4$}     at 6 2.2   %
\put {$^3$}     at 5 2.2   %
\put {$^2$}     at 4 2.2   %
\put {$^1$}     at 3 2.2   % label nodes with roots below
\linethickness=0.75pt                          % sets line thickness
\putrule from 4.03 2.045 to 4.97 2.045       % puts solid lines between nodes
\putrule from 4.03 1.955 to 4.97 1.955       %
\putrule from 3.05 2 to 3.95 2              %
\putrule from 5.05 2 to 5.95 2              %
%\putrule from 6.05 2 to 6.95 2              %
}
\endpicture
\qquad\qquad
%********************************************************************
% Dynkin diagram of type H_3
%********************************************************************
\beginpicture
\setcoordinatesystem units <1cm,1cm>        % sets scale
\setplotarea x from 2 to 5.2, y from 1.5 to 2.8  % sets plot size up
\put{$H_3$} at 2 2.5
{\scriptsize
\multiput {$\circ$} at 3   2 *2 1 0 /      %  puts nodes in
\put {$^3$}     at 5 2.2   %
\put {$^2$}     at 4 2.2   %
\put {$^5$}     at 3.5 2.1 %
\put {$^1$}     at 3 2.2   % label nodes with roots below
\linethickness=0.5pt                          % sets line thickness
%\putrule from 3.03 2.045 to 3.97 2.045       % puts solid lines between nodes
%\putrule from 3.03 1.955 to 3.97 1.955       %
\putrule from 3.05 2 to 3.95 2              %
\putrule from 4.05 2 to 4.95 2              %
}
% end of dynkin diagram of type H_3
\endpicture
$$
$$
%********************************************************************
% Dynkin diagram of type H_4
%********************************************************************
\beginpicture
\setcoordinatesystem units <1cm,1cm>        % sets scale
\setplotarea x from 3.5 to 7.2, y from 1.5 to 2.8  % sets plot size up
\put{$H_4$} at 2 2.5
{\scriptsize
\multiput {$\circ$} at  3  2 *3 1 0 /      %  puts nodes in    
\put {$^4$}     at 6 2.2   %
\put {$^3$}     at 5 2.2   %
\put {$^2$}     at 4 2.2   %
\put {$^5$}   at 3.5 2.1 %
\put {$^1$}     at 3 2.2   % label nodes with roots below
\linethickness=0.5pt                          % sets line thickness
\putrule from 3.05 2 to 3.95 2              %
\putrule from 4.05 2 to 4.95 2              %
\putrule from 5.05 2 to 5.95 2
}
% end of dynkin diagram of type H_4
\endpicture
\qquad\qquad
%********************************************************************
% Dynkin diagram of type  I_2(r)
%********************************************************************
\beginpicture
\setcoordinatesystem units <1cm,1cm>       % sets scale
\setplotarea x from 1.8 to 4.2, y from 1.5 to 2.8 % sets plot size up
\put{$I_2(m)$} at 2 2.5
{\scriptsize
\multiput {$\circ$} at 3   2 *1 1 0 /      % puts nodes in
\put {$^1$}     at 3 2.2   %
\put {$^m$}     at 3.5 2.1   %
\put {$^2$}     at 4 2.2   % label nodes with roots above
\linethickness=0.5pt                      % sets line thickness
%\putrule from 3.05 0 to 3.95 0              %
%\putrule from 3.03 2.045 to 3.97 2.045  %
%\putrule from 3.03 1.955 to 3.97 1.955  %
\putrule from 3.05 2 to 3.95 2         %   puts solid lines between nodes
}
% end of dynkin diagram of type I_2(r)
\endpicture
$$
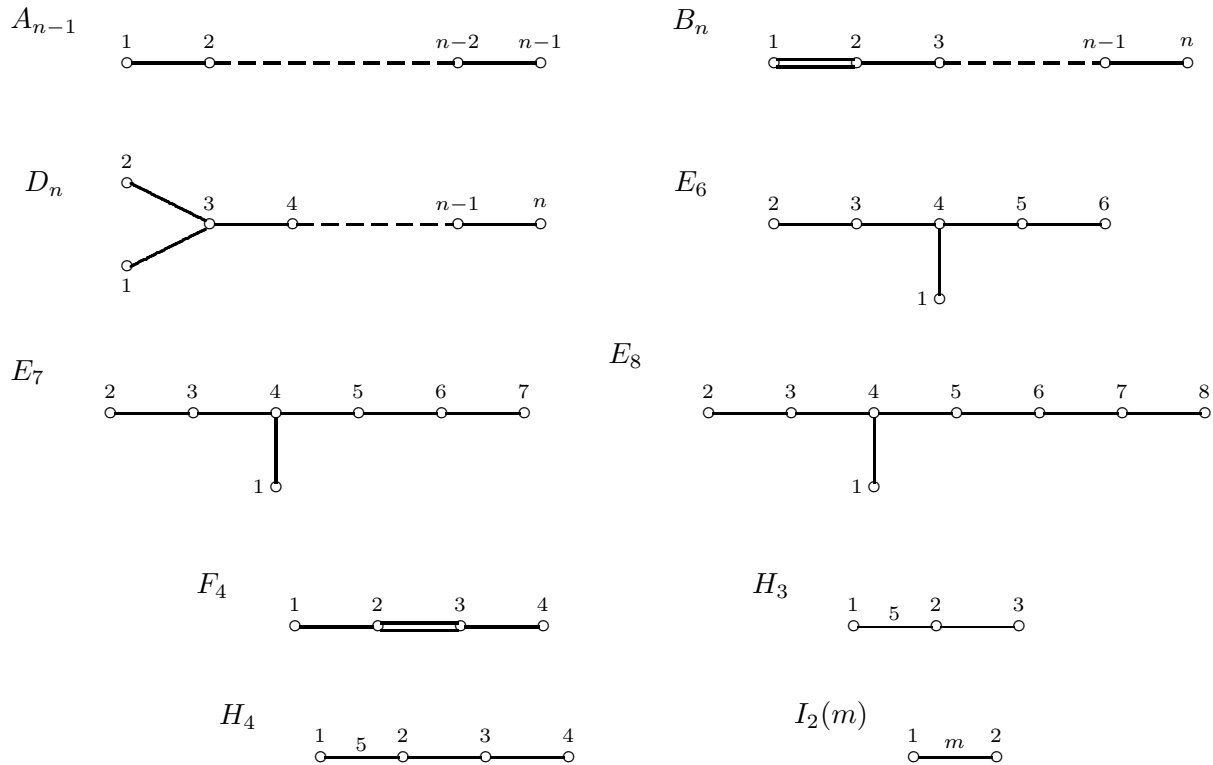
\centerline{{\bf Figure 1.}  Coxeter-Dynkin diagrams for real reflection groups.}

\vfill\eject

\subsection 2B.\ \ Complex reflection groups

The irreducible complex reflection groups were classified
by Shephard and Todd [ST].  There is one infinite family denoted
$G(r,p,n)$ and a list of exceptional complex reflection groups
denoted $G_4,\ldots, G_{35}$.  In this subsection, we classify the 
graded Hecke algebras for the groups $G(r,p,n)$.

Let $r$, $p$ and $n$ be positive integers with $p$ dividing $r$ and
let $\xi=e^{2\pi i/r}$. 
Let $S_n$ be the symmetric group of $n\times n$ matrices and let 
$$\xi_j = \diag(1,1,\ldots, 1, \xi, 1,\ldots,1),$$
where $\xi$ appears in the $j$th entry. 
Then
$$ 
G(r,p,n) = \{
\xi_1^{\lambda_1}\cdots \xi_n^{\lambda_n}w\ |\ 
w\in S_n, \ 0 \leq \lambda_i \le r-1,\ \lambda_1+\cdots+ \lambda_n=0 \bmod p\}.$$
For $\lambda = (\lambda_1,\ldots, \lambda_n)\in (\ZZ/r\ZZ)^n$, let
$\xi^\lambda = \xi_1^{\lambda_1}\cdots\xi_n^{\lambda_n}$.  Then the
multiplication in $G(r,p,n)$ is described by the relations
$$\xi^\lambda\xi^\mu=\xi^{\lambda+\mu}
\qquad\hbox{and}\qquad
w\xi^\lambda=\xi^{w\lambda}w,
\qquad\hbox{for $\lambda,\mu\in (\ZZ/r\ZZ)^n$, $w\in S_n$,}$$
where $S_n$ acts on $(\ZZ/r\ZZ)^n$ by permuting the factors.
Let $v_i$ be the column vector with $1$ in the $i^{\rm th}$ entry and all other 
entries $0$.
The group $G(r,p,n)$ acts on $V:= \CC^n$ with orthonormal basis
$\{ v_1, \ldots, v_n \}$ as a complex reflection group.

Every real reflection group is a complex reflection group and several
of these are special cases of the groups $G(r,p,n)$.  In particular,
\itemitem{(a)} $G(1,1,n)$ is the symmetric group $S_n$,
\itemitem{(b)} $G(2,1,n)$ is the Weyl group $WB_n$ of type $B_n$,
\itemitem{(c)} $G(2,2,n)$ is the Weyl group $WD_n$ of type $D_n$, and 
\itemitem{(d)} $G(r,r,2)$ is the dihedral group $I_2(r)$ of order $2r$.

The reflections in $G(r,p,n)$ are 
$$\matrix{
\xi^{kp}_i, &\qquad\qquad &1\le i\le n,\enspace  0\le k\le (r/p)-1,
\enspace \hbox{and} \cr
\cr
\xi^{k}_i\xi_j^{-k}(i,j), &\qquad\qquad &1\le i<j\le n,\enspace  0\le k\le r-1, \cr
}$$
where $(i,j)$ is the transposition in $S_n$ that switches $i$ and $j$.

\smallskip\noindent
{\bf Conjugacy in $G(r,p,n)$.}  Each element of $G(r,p,n)$ is
conjugate by elements of $S_n$ to a disjoint product of cycles of the form
$$\xi_i^{\lambda_i}\cdots \xi_k^{\lambda_k}(i,i+1,\ldots,k).$$
By conjugating this cycle by 
$\xi_i^{-c}\xi_{i+1}^{\lambda_i}\xi_{i+2}^{\lambda_i+\lambda_{i+1}}\cdots
\xi_k^{\lambda_i+\cdots + \lambda_{k-1}}\in G(r,r,n)$, we have
$$\xi_i^{-c}\xi_k^{c+\lambda_i+\cdots+\lambda_k}(i,\ldots,k),
\qquad\hbox{where}\ \ c=(k-i)\lambda_i+(k-i-1)\lambda_{i+1}+\cdots+\lambda_{k-1}.$$
If $i_1,i_2,\ldots, i_\ell$ denote the minimal indices of the 
cycles and $c_1,\ldots,c_\ell$ are the numbers $c$ for the various
cycles, 
then after conjugating by
$\xi_{i_1}^{c_1}\cdots \xi_{i_{\ell-1}}^{c_{\ell-1}}
\xi_{i_{\ell}}^{-(c_1+\cdots+c_{\ell-1})}\in G(r,r,n)$, 
each cycle becomes
$$\xi_k^{\lambda_i+\cdots+\lambda_k}(i,\ldots,k)
\qquad\hbox{except the last, which is}\qquad
\xi_{i_\ell}^{-a}\xi_n^b(i_{\ell},\ldots,n),$$
where $a=c_1+\cdots+c_\ell$ and $b=a+\lambda_{i_\ell}+\cdots+\lambda_n$.
If $k=n-i_\ell+1$ is the length of the last cycle, 
then conjugating the last cycle by
$\xi_{i_\ell}^{k-1}\xi_{i_{\ell}+1}^{-1}\cdots \xi_n^{-1}\in G(r,r,n)$
gives
$$\xi_{i_\ell}^{-a+k}\xi_n^{b-k}(i_{\ell},\ldots,n).$$
If we conjugate the last cycle by $\xi_{i_\ell}^p\in G(r,p,n)$,
we have 
$$\xi_{i_\ell}^{-a+p}\xi_n^{b-p}(i_{\ell},\ldots,n).$$
In summary, any element $g$ of $G(r,p,n)$ is conjugate 
to a product of disjoint
cycles where each cycle is of the form
$$\xi_k^a(i,i+1,\ldots,k), \qquad 0\le a\le r-1,
\global\advance\resultno by 1\eqno{(\the\sectno.\the\resultno {\rm a})}$$
except possibly the last cycle, which is of the form 
$$\xi_{i_\ell}^a\xi_n^b(i_\ell,i_\ell+1,\ldots,n),
\qquad\hbox{with $0\le a\le \gcd(p,k)-1$,}
\eqno{(\the\sectno.\the\resultno {\rm b})}$$
where $k=n-i_{\ell}+1$ is the length of the last cycle.

\medskip\noindent
{\bf Centralizers in $G(r,p,n)$.}  Let
$Z_{G(r,p,n)}(g)=\{ h\in G(r,p,n)\ |\ hg=gh\}$ denote the centralizer of
$g\in G(r,p,n)$.  Since $G(r,p,n)$ is a subgroup of $G(r,1,n)$, 
$$Z_{G(r,p,n)}(g) = Z_{G(r,1,n)}(g)\cap G(r,p,n),$$
for any element $g\in G(r,p,n)$.  Suppose that
$g$ is an element of $G(r,1,n)$ which is a product of disjoint
cycles of the form $\xi_k^a(i,\ldots, k)$ and that $h\in G(r,1,n)$ 
commutes with $g$. Conjugation by $h$ effects some
combination of the following operations on the cycles of $g$:
\item{(a)} permuting cycles of the same type,
$\xi_k^a(i,\ldots,k)$
and $\xi_m^b(j,\ldots, m)$ with $b=a$ and $k-i=m-j$,
\item{(b)} conjugating a single cycle 
$\xi_k^a(i,\ldots,k)$ by powers of itself, and
\item{(c)} conjugating a single cycle 
$\xi_k^a(i,\ldots,k)$ by $\xi_i^b\cdots \xi_k^b$,\ \  for any 
$0\le b\le r-1$.
\smallskip\noindent
Furthermore, the elements of $G(r,1,n)$ which commute with $g$ are determined
by how they ``rearrange'' the cycles of $g$ and 
a count (see [Mac, p.~170]) of the number of such operations shows that
if $g\in G(r,1,n)$ and $m_{a,k}$ is the number of cycles of type
$\xi_{i+k}^a(i,i+1,\ldots,i+k)$ for $g$, then 
$$\Card(Z_{G(r,1,n)}(g)) 
= \prod_{a,k} (m_{a,k}!\cdot k^{m_{a,k}}\cdot r^{m_{a,k}})\,.
\formula$$

\medskip\noindent
{\bf Determining the graded Hecke algebras for $G(r,p,n)$.}
It follows from Lemma 1.8a that 
if $g=\xi_i^{a+b}\xi_k^{-a}(i,\ldots,k)$, then $(V^g)^\perp$ has basis
$$
\{v_k-v_{k-1},v_{k-1}-v_{k-2},\ldots,v_{i+1}-\xi^a v_i\} \quad \hbox{if $b=0$,} 
\qquad\hbox{and} \qquad 
\{v_i,\ldots, v_k\} \quad\hbox{if $b\ne 0$.} 
$$
Thus, if $g\in G(r,p,n)$ and $\codim(V^g)=2$, then $g$ is conjugate to
one of the following elements:
$$\matrix{
b= \xi_1^{a}\xi_3^{-a}(1,2,3), \hfill
&\qquad &\hbox{$0\le a\le \gcd(p,3)-1$,} \hfill \cr
\cr
c=\xi_1^{a+\ell}\xi_2^{-a}(1,2),\hfill 
&\qquad &\hbox{$\ell\neq 0$ (so $r\ne 1$),} \hfill \cr
\cr
d=\xi_1^{\ell_1}\xi_2^{\ell_2},\hfill 
&&\hbox{$\ell_1 \neq 0$, $\ell_2\ne 0$ (so $r\ne 1$)}, \hfill\cr
\cr
e=(1,2)\xi_3^\ell,\hfill && \ell \neq 0, \hfill\cr 
\cr
f=(1,2)\xi_3^a\xi_4^{-a}(3,4). \cr
}
$$
It is interesting to note that these elements are also
representatives of the conjugacy classes of elements in $G(r,p,n)$
which can be written as a product of two reflections.

We determine conditions on the above elements and on $r$, $p$, and $n$
to give nontrivial graded Hecke algebras:
\item{(z)}  The center of $G(r,p,n)$ is 
$$Z(G(r,p,n)) = \{ \xi_1^\ell\cdots\xi_n^\ell\ |\ n\ell = 0 \bmod p\}.$$
Since $\xi_1^p\cdots\xi_n^p\in Z(G(r,p,n))$, it follows that $p=r$ or
$p=r/2$ whenever $Z(G(r,p,n))\subseteq \{\pm 1\} 
=\{\xi_1^0\cdots \xi_n^0,\xi_1^{r/2}\cdots\xi_n^{r/2}\}$.
\item{(b1)} If $n\ge 4$, the element $\xi_1\xi_2\xi_3\xi_4^{-3}\in Z_G(b)$
and has determinant $\xi^2$ on $(V^b)^\perp=\hbox{span-}\{v_3-v_2,v_2-\xi^av_1\}$.
\item{(b2)} If $n=3$ and $p=0$ mod $3$, the element 
$\xi_1^{p/3}\xi_2^{p/3}\xi_3^{p/3}\in Z_G(b)$ and has determinant $\xi^{2p/3}$
on $(V^b)^\perp$.
\item{(c1)} If $n\ge 3$, the element $\xi_1\xi_2\xi_3^{-2}\in Z_G(c)$
and has determinant $\xi^2$ on $(V^c)^\perp=\hbox{span-}\{v_1,v_2\}$.
\item{(c2)} If $n=2$, $p=r/2$ and $p$ is odd, the element
$\xi_1^{p/4}\xi_2^{p/4}\in Z_G(c)$ 
and has determinant $\xi^{r/2}$ on $(V^c)^\perp$.
\item{(d1)} If $n\ge 3$, the element $\xi_1\xi_3^{-1}\in Z_G(d)$
and has determinant $\xi$ on $(V^d)^\perp=\hbox{span-}\{v_1,v_2\}$.
\item{(d2)} If $p=r/2$, the element $\xi_1^{r/2}\in Z_G(d)$
and has determinant $\xi^{r/2}$ on $(V^d)^\perp$.
\item{(ef)} The elements $e$ and $f$ have order $2$.
\smallskip\noindent
Thus, it follows from Corollary 1.11, Theorem 1.9b, and Lemma 2.2b that 
if $A$ is a graded Hecke algebra for $G(r,p,n)$, then
$$\matrix{
a_b=0 \hfill &\quad &\hbox{unless} &\quad &\hbox{(i) $r=1$, or} \hfill\cr
&&&&\hbox{(ii) $r=2$, or} \hfill\cr
&&&&\hbox{(iii) $n=3$ and $p\ne 0 \bmod 3$,} \hfill\cr
\cr
a_c=0 \hfill &&\hbox{unless} &&\hbox{(i) $r=2$ and $p=1$, or} \hfill\cr
&&&&\hbox{(ii) $n=2$ and $p=r/2$,} \hfill\cr
\cr
a_d=0 \hfill &&\hbox{unless} &&\hbox{$p=r$, $n=2$ and $p\ne 0\bmod 2$,} \hfill\cr
\cr
a_e=0 \hfill &&\hbox{always,} &&\hbox{and} \hfill\cr
\cr
a_f=0 \hfill &&\hbox{always.} \hfill\cr
}$$
In the remaining cases, one uses the description of $Z_G(g)$ given 
just before (2.5) to check that all elements of $Z_G(g)$ have 
determinant $1$ on $(V^g)^\perp$.
Note that $n=3$ and $p\ne 0 \bmod 3$ imply that $a_b=0$  
for the elements $b= \xi_1^{a}\xi_3^{-a}(1,2,3)$.

We arrive at the following enumeration of the nontrivial
graded Hecke algebras for complex reflection groups. (The tensor
product algebra $S(V)\otimes \CC G$ always exists and corresponds
to the case when all of the skew symmetric forms $a_g$ are zero).  
The table below gives representatives of the conjugacy classes
of $g\in G$ that may have $a_g\ne 0$ for some graded Hecke algebra $A$.

$$
\vbox{
\offinterlineskip
\hrule
\halign{& \vrule# & 
                    \quad\  $#$ \hfil\quad
        & \vrule# &  \quad\hfil $#$ \hfil\quad 
        & \vrule#  \cr
 & \strutH \ \ \quad \ {\rm Group}\ && {\rm Representative}\  g & \cr
 & \strutL  &&{\rm with}\ a_g \ne 0 & \cr
\noalign{\hrule}
&  \strutB  G(1,1,n)=S_n   && (1,2,3) & \cr
& \strutB G(2,1,n)=WB_n,\ \ n\ge 3  && \xi_1(1,2),\ (1,2,3) &\cr
&  \strutB  G(2,2,n)=WD_n,\ \ n\ge 3 && (1,2,3) &\cr
&\strutB G(r,r,2)=I_2(r) &&  \xi_1^k \xi_2^{r-k}, \ 0<k<r/2  &\cr
&\strutB G(r,r/2,2), \ \ r/2\ {\rm odd}  && \xi_2^{r/2}(1,2) &\cr
& \strutB G(r,r,3), \ \ r\ne 0 \bmod 3 && (1,2,3) &\cr
& \strutB G(r,r/2,3),\  \ r/2 \ne 0 \bmod 3,\  r\ne 2 && (1,2,3) &\cr
}
\hrule}
$$
\centerline{{\bf Table 2.}  Graded Hecke algebras for the groups $G(r,p,n)$.}

\subsection 2C.  Exceptional complex reflection groups

The irreducible
exceptional complex reflection groups $G$ are denoted $G_4,\ldots, G_{35}$ 
in the classification of Shephard and Todd.  From Table VII in 
[ST], one sees that the center of $G$ is $\pm 1$ only in the cases
$G_4$, $G_{12}$, $G_{24}$ and $G_{33}$.  
By Schur's lemma, the center of an irreducible complex reflection group 
consists of multiples of the identity.  Thus, by
Corollary 1.11, the only exceptional complex reflection groups
that could have a nontrivial graded Hecke algebra
(i.e., with some $a_g\ne 0$) 
are $G_4$, $G_{12}$, $G_{24}$ and $G_{33}$
(we exclude the real groups).
We determine the graded Hecke algebras for these groups using
Theorem 1.9b and Lemma 2.2.

\medskip
The rank $2$ group $G_4$ has order 24 and seven conjugacy classes.
The following data concerning these conjugacy classes are
obtained from the computer software GAP [S+]
using the package CHEVIE [G+].
In the following table, $\omega$ is a primitive cube root of unity
and $C(g)$ denotes the conjugacy class of $g$.  
$$
\vbox{
\offinterlineskip
\hrule
\halign{& \vrule# & 
     \strutA %homemade strut
                    \quad\hfil $#$ \quad
        & \vrule# & \quad\hfil $#$ \quad
        & \vrule# & \quad\hfil $#$ \quad
        & \vrule# & \quad\hfil $#$ \quad
        & \vrule# & \quad\hfil $#$ \quad
        & \vrule# & \quad\hfil $#$ \quad
        & \vrule# & \quad\hfil $#$ \quad
        & \vrule# & \quad\hfil $#$ \ \ 
         & \vrule# \cr
&& \multispan{14}\hfil 
     {${\rm Conjugacy\ class\ representatives\ for}\ G_{4}$} 
   \qquad\qquad\qquad\hfil &\cr
\noalign{\hrule}
& {\rm Order}(g) && 1 && 4 && 3 && 6 && 6 && 3 && 2 & \cr
 \noalign{\hrule}
& \det(g) && 1 && 1 && \omega && \omega && \omega^{-1}&& \omega^{-1} && 1  & \cr
 \noalign{\hrule}
&  \hbox{$\vert$}\cC(g)\hbox{$\vert$} 
     && 1 && 6 && 4 && 4 && 4 && 4 && 1 &\cr
\noalign{\hrule}
&\hbox{$\vert$}Z_G(g)\hbox{$\vert$}  
&& 24 && 4 && 6 && 6 && 6 && 6 && 24 &\cr
}
\hrule}
$$
The elements with determinant $1$ and order more than $2$ 
in $G_4$ all have order $4$.  If $g$ is an element of order $4$,
then $|Z_G(g)| = 4$ and
every element of $Z_G(g)$ has determinant $1$ since
$Z_G(g)$ is generated by $g$.
Hence, by Theorem 1.9b and Lemma 2.2, $a_g$ can be nonzero for
a graded Hecke algebra for $G_4$ exactly when $g$ has order $4$.
Thus, the dimension of the space of parameters for graded Hecke algebras 
of $G_{4}$ is 1.

\medskip
The rank $2$ group $G_{12}$ has order $48$.  The computer software GAP
provides the following information about the conjugacy classes
of $G_{12}$.
$$
\vbox{
\offinterlineskip
\hrule
\halign{& \vrule# & 
     \strutA %homemade strut
                    \quad\hfil $#$ \quad 
        & \vrule# & \quad\hfil $#$ \quad 
        & \vrule# & \quad\hfil $#$ \quad 
        & \vrule# & \quad\hfil $#$ \quad 
        & \vrule# & \ \quad\hfil $#$ \quad 
        & \vrule# & \quad\hfil $#$ \quad 
        & \vrule# & \quad\hfil $#$ \quad 
        & \vrule# & \ \quad\hfil $#$ \quad 
        & \vrule# & \quad\hfil $#$ \quad 
        & \vrule# \cr
% \strut \ \ #  \cr
%height2pt&\omit &\omit & \cr
&& \multispan{16}\hfil 
      Conjugacy class representatives for $G_{12}$ 
   \qquad\qquad\qquad\hfil &\cr
\noalign{\hrule}
& \hbox{Order($g$)} && 1 && 2 && 8 && 6 && 8 && 2 && 3 && 4 & \cr
 \noalign{\hrule}
& \hbox{det}(g) && 1 && -1 && -1 && 1 && -1 && 1 && 1 && 1 & \cr
 \noalign{\hrule}
&  \hbox{$\vert$}\cC(g)\hbox{$\vert$} 
     && 1 && 12 && 6 && 8 && 6 && 1 && 8 && 6 &\cr
\noalign{\hrule}
&\hbox{$\vert$}Z_G(g)\hbox{$\vert$}  
%= \hbox{$\vert$}$G$\hbox{$\vert$}/\hbox{$\vert$}$\cC(G)$\hbox{$\vert$}
&& 48 && 4 && 8 && 6 && 8 && 48 && 6 && 8 &\cr
}
\hrule}
$$
If $g$ is an element in $G_{12}$ with order more than $2$ and determinant $1$,
then $g$ has order $3$, $4$, or $6$.
Let $h$ be any element of order $8$.
Then $h$ has determinant $-1$ and
commutes with $h^2$ of order $4$.
Hence, by Theorem 1.9b, if $g$ has order $4$,
then $a_g=0$. 
Let $g_6$ be a representative from the class of elements
of order $6$.  Since $|Z_G(g_6)| = 6$,
$Z_G(g_6)$ is generated by $g_6$ and hence
every element of $Z_G(g_6)$ has determinant 1. 
We can choose $g_6^2$ as a representative
for the conjugacy class of elements of order $3$.
As $g_6$ and $g_6^2$ commute,
$\langle g_6 \rangle \subset Z_G(g_6^2)$.
But $|\langle g_6 \rangle| = 6 = |Z_G(g_3)|$,
so $Z_G(g_6^2)$ is generated by $g_6$
and every element of $Z_G(g_3)$ has determinant $1$.
Thus, $a_g$ can be nonzero for a graded Hecke algebra $A$ for
$G_{12}$ exactly when $g$ has order $3$ or $6$.
Thus, the dimension of the space of parameters of graded Hecke algebras for 
$G_{12}$ is $2$.

\medskip
The rank $3$ group $G_{24}$ has order $336$.  Note that
$-1\in G_{24}$ since $Z(G)=\{\pm 1\}$.
Up to $G$-orbits, there are two codimension 2 subspaces, $L$ and
$M$, that are equal to $V^g$ for some $g\in G_{24}$
(see [OT, App.~C, Table C.5]).
Furthermore, ${\rm Stab}(L) \cong A_2$ and ${\rm Stab}(M) \cong B_2$.  
We need only
consider elements of order 3 in ${\rm Stab}(L) \cong A_2$ and 
of order 4 in ${\rm Stab}(M)\cong B_2$ (as the rest have order
$1$ or $2$).
In $G_{24}$, there is only one conjugacy class of elements
of order $3$ and only one conjugacy class of elements of order $4$
and determinant $1$.  The table below (obtained using GAP) records
information about these classes.
\bigskip
$$
\vbox{
\offinterlineskip
\hrule
\halign{& \vrule# & 
     \strutA %homemade strut
                    \quad\hfil $#$ \quad 
        & \vrule# & \quad\hfil $#$ \quad 
        & \vrule# & \quad\hfil $#$ \quad 
        & \vrule# \cr
%  \strut \ \ #  \cr
%\height2pt&\omit &\omit & \cr
& \multispan{5}
\hfil 
\strutA Certain classes of 
      $G_{24}$ 
  \hfil 

 &  \cr
\noalign{\hrule}
& \hbox{Order($g$)} && 3 && 4 & \cr
 \noalign{\hrule}
& \hbox{det}(g) && 1 && 1 & \cr
 \noalign{\hrule}
&  \hbox{$\vert$}\cC(g)\hbox{$\vert$} 
     && 56 && 42 &\cr
\noalign{\hrule}
&\hbox{$\vert$}Z_G(g)\hbox{$\vert$}  
%= \hbox{$\vert$}$G$\hbox{$\vert$}/\hbox{$\vert$}$\cC(G)$\hbox{$\vert$}
&& 6 && 8 &\cr
}
\hrule}
$$
\bigskip\noindent
If $g$ has order $3$,
$Z_G(g)$ must 
contain $1,g,g^2$, and $-1$,
and hence $Z_G(g)$ is generated by these elements
since $|Z_G(g)|=6$.
Thus all elements of $Z_G(g)$
have determinant $1$ on $(V^g)^\perp$.  
If $g$ has order 4 and determinant $1$,
then $Z_G(g)$ must contain $1,g,g^2,g^3$, and $-1$,
elements which all have determinant $1$ on $(V^g)^\perp$.
Since $|Z_G(g)|=8$, these elements generate $Z_G(g)$
and so every element of $Z_G(g)$ has determinant $1$ on $(V^g)^\perp$.
Hence, $a_g$ can be nonzero for a graded Hecke algebra of $G_{24}$
exactly when $g$ has order $3$ or $g$ has order $4$ and determinant $1$.
Thus, the dimension of the space of parameters for
graded Hecke algebras for $G_{24}$ is $2$.

\medskip
The group $G_{33}$ is the only exceptional complex
reflection group of rank $5$. It has order $72\cdot 6!$ and 
degrees $4,6,10,12, 18$.  
There are $45$ reflecting hyperplanes and the corresponding
reflections all have order $2$.
Up to $G$-orbits, there are two codimension
$2$ subspaces, $L$ and $M$,
that are equal to $V^g$ for some $g\in G_{33}$
(see [OT, App.~C, Table C.14]).
Furthermore,
${\rm Stab}(L) \cong A_1 \times A_1$ and ${\rm Stab}(M) \cong A_2$.
We need not consider the case where $V^g = L$
since then $g$ has order $2$ and 
hence $a_g=0$ for any graded Hecke algebra
by Proposition~2.2b.

We use a presentation for $G_{33}$ in six coordinates from [ST]:
Let $V:= \{ (x_1,x_2,x_3,x_4,x_5,x_6)\ | \ x_i\in \CC\}$
and consider the group generated by order $2$ reflections about the hyperplanes
$H_1 = \{x_2 - x_3=0\}$,
$H_2 = \{x_3 - x_4=0\}$,
$H_3 = \{x_1 - x_2=0\}$,
$H_4 = \{x_1 - \omega x_2=0\}$,
$H_5 = \{x_1 + x_2 + x_3 + x_4 + x_5 + x_6 = 0\}$,
where $\omega$ is a primitive cube root of unity.
The fixed point space of this (reducible) group is
$Y = H_1\cap\cdots\cap H_5=\{(0,0,0,0,x,-x)\ |\ x\in \CC\}$,
and $G_{33}$ is just the restriction to $Y^\perp$.
Let $s_i$ be the order $2$ reflection about $H_i$.
Let $g = s_1 s_3$. Then $V^g$ = $H_1 \cap H_3$ and 
 ${\rm Stab}(V^g) \cong A_2$.
Let 
$ h = (s_1 s_3 s_4)^2  $, 
the diagonal matrix with diagonal 
$\{\omega,\omega,\omega,1,1,1\}$.
Then $h$ acts as $\omega$ times the identity 
on $(V^g)^\perp$ as $(V^g)^\perp \subseteq
\hbox{$\CC$-span}\{x_1,x_2, x_3\}$.
Hence, $h$ commutes with $g$. 
But $(V^g)^\perp$ has dimension $2$ and
$h$ has determinant $\omega^2 \neq 1$ on 
$(V^g)^\perp$.  Thus, by Theorem 1.9b and Lemma 2.2, $a_g=0$
for any graded Hecke algebra.  The same argument applied to $Y^\perp$ shows
that $G_{33}$ has no nontrivial graded Hecke algebras. 
In summary, the dimension of the space of parameters for
graded Hecke algebras for $G_{33}$ is zero.

\medskip
$$
\vbox{
\offinterlineskip
\hrule
\halign{& \vrule# & 
                    \quad\  $#$ \hfil\quad
        & \vrule# &  \quad\hfil $#$ \hfil\quad 
        & \vrule#  \cr
 & \strutA \ \ \quad \ {\rm Group}\ && \hbox{$g$ with $a_g\ne 0$} & \cr
\noalign{\hrule}
&  \strutB G_4     && {\rm Order}(g)=4 & \cr
&  \strutB G_{12}  && {\rm Order}(g_1)=3\ \hbox{and}\  {\rm Order}(g_2)=6 &\cr
&  \strutB G_{24}  && {\rm Order}(g_1)=3\ \hbox{and}\ 
{\rm Order}(g_2)=4, \det(g_2)=1 &\cr
}
\hrule}
$$
\centerline{{\bf Table 3.}  Graded Hecke algebras for 
exceptional complex reflection groups.}

\section 3. \ The graded Hecke algebras $H_{\rm gr}$

In [Lu], Lusztig gives a definition of graded Hecke algebras
for real reflection groups which is different from
the definition in Section 1,
which applies to more general groups.
It is not obvious that Lusztig's algebras are
examples of the graded Hecke algebras defined
in Section 1.  In this section, we show explicitly
how the definition of Section 1 includes Lusztig's algebras.

Let $W$ be a finite real reflection group acting on $V$
and let $R$ be the root system of $W$. Let
$\alpha_1,\ldots, \alpha_n$ be a choice of simple roots
in $R$ and let $s_1,\ldots, s_n$ be the corresponding simple 
reflections in $W$.
Let $s_\alpha$ be the reflection in the root $\alpha$ so that,
for $v\in V$,
$$s_\alpha v = v - \langle v,\alpha^\vee\rangle \alpha,
\qquad\hbox{where $\alpha^\vee = 2\alpha/\langle\alpha,\alpha\rangle$.}$$ 
Let $R^+=\{\alpha>0\}$ denote the set of positive roots in $R$.

Let $k_\alpha$ be fixed complex numbers indexed by the
roots $\alpha\in R$ satisfying
$$k_{w\alpha} = k_\alpha,
\qquad\hbox{for all $w\in W$, $\alpha\in R$.}\formula$$
This amounts to a choice of either one or two ``parameters'', depending
on whether all roots in $R$ are the same length or not.
As in Section 1, let
$\CC W= \hbox{$\CC$-span}\{t_g \ |\ g \in W\}$,
with $t_g t_h = t_{gh}$,
and let $S(V)$ be the symmetric algebra of $V$.  Lusztig [Lu] defines the ``graded Hecke algebra'' with parameters
$\{k_\alpha\}$ to be the unique algebra structure $H_{\rm gr}$ on the
vector space $S(V)\otimes \CC W$ such that 
\smallskip\noindent
\global\advance\resultno by 1
\itemitem{(\the\sectno.\the\resultno~a)}  
$S(V)=S(V)\otimes 1$ is a subalgebra of $H_{\rm gr}$,
\smallskip\noindent
\itemitem{(\the\sectno.\the\resultno~b)}  
$\CC W= 1\otimes \CC W$ is a subalgebra of $H_{\rm gr}$, and
\smallskip\noindent
\itemitem{(\the\sectno.\the\resultno~c)} $t_{s_i} v = (s_i v)t_{s_i} - k_{\alpha_i}
\langle v,\alpha_i^\vee\rangle,$
for all $v\in V$ and simple reflections $s_i$ in the simple roots $\alpha_i$.  

\smallskip\noindent
We shall show that every algebra $H_{\rm gr}$ as defined by (3.2~a-c) is 
a graded Hecke algebra $A$ for a specific set of skew symmetric bilinear 
forms $a_g$.

Let $k_\alpha\in \CC$ as in (3.1).  Use the notation
$$h = \hbox{${1}\over {2}$} \sum_{\alpha>0} k_\alpha\alpha^\vee t_{s_\alpha},
\qquad\hbox{so that}\qquad
\langle v,h\rangle = \hbox{$1\over2$}
\sum_{\alpha>0} k_\alpha\langle v,\alpha^\vee\rangle
t_{s_\alpha}
\formula
$$
for $v\in V$.  The element $h$ should be viewed
as an element of $V\otimes \CC W$, and 
$\langle v,h\rangle\in \CC W$.  With this notation, 
let $A$ be the algebra (as in Section 1) generated by
$V$ and $\CC W$ with relations
$$t_g v = (gv)t_g
\qquad\hbox{and}\qquad
[v,w]=-[\langle v,h\rangle,\langle w,h\rangle],
\qquad\hbox{for $v,w\in V$, $g\in W$.} 
\formula
$$
Note that $A$ is defined by the bilinear forms
$$
a_g(v,w) = 
\hbox{$1\over 4$}
\sum_{
\hbox{$\alpha ,\, \beta \, > \,0 \atop
g\, =\, s_\alpha s_\beta$}}
k_\alpha k_\beta
\left(
\langle v, \beta^\vee \rangle\langle w, \alpha^\vee \rangle
-\langle v, \alpha^\vee \rangle\langle w, \beta^\vee \rangle
\right).
$$

The following theorem shows that the algebra $A$ satisfies the defining
conditions (3.2 a-c) of the algebra $H_{\rm gr}$.

\thm  Let $W$ be a finite real reflection group and let $A$ be the
algebra defined by (3.4).
\item{(a)}  As vector spaces, $A\cong S(V)\otimes \CC W$
(and hence, $A$ is a graded Hecke algebra).
\item{(b)}  If $\tilde v=v-\langle v,h\rangle$ for $v\in V$,
then
$$[\tilde v,\tilde w]=0\qquad\hbox{and}\qquad
t_{s_i}\tilde v = (\widetilde{s_i v})t_{s_i} - k_{\alpha_i}
\langle v,\alpha_i^\vee\rangle,$$
for all $v,w\in V$ and simple reflections $s_i$ in $W$.
\pf
First note that if $u,v\in V$ then
$$[u,\langle v,h\rangle]=
\hbox{$1\over 2$} \sum_{\alpha>0} k_\alpha \langle v,\alpha^\vee\rangle
\langle u,\alpha^\vee\rangle \alpha t_{s_\alpha}
=[v,\langle u,h\rangle].\eqno{(*)}$$
Thus, for $u,v,w\in V$,
$$\eqalign{
[u,[v,w]]+&[w,[u,v]]+[v,[w,u]] \cr
&=[u, [\langle w,h\rangle,\langle v,h\rangle]]
+[w, [\langle v,h\rangle,\langle u,h\rangle]]
+[v, [\langle u,h\rangle,\langle w,h\rangle]] 
\cr
&=[[u, \langle w,h\rangle],\langle v,h\rangle]
+[\langle w,h\rangle,[u,\langle v,h\rangle]]
+[[w, \langle v,h\rangle],\langle u,h\rangle] \cr
&\qquad\qquad
+[\langle v,h\rangle,[w,\langle u,h\rangle]] 
+[[v, \langle u,h\rangle],\langle w,h\rangle] 
+[\langle u,h\rangle],[v,\langle w,h\rangle]]  
\cr
&=[[w, \langle u,h\rangle],\langle v,h\rangle]
+[\langle w,h\rangle,[v,\langle u,h\rangle]]
+[[v, \langle w,h\rangle],\langle u,h\rangle] \cr
&\qquad\qquad
+[\langle v,h\rangle,[w,\langle u,h\rangle]] 
+[[v, \langle u,h\rangle],\langle w,h\rangle] 
+[\langle u,h\rangle,[v,\langle w,h\rangle]]  
\cr
&=0. \cr}
\formula$$
For $v\in V$, $h\in W$, and $s_i$ a simple reflection,
$$\eqalign{
t_{s_i} \langle v, h \rangle t_{s_i}
&={1\over2}\sum_{\alpha>0} k_\alpha\langle v,\alpha^\vee\rangle
t_{s_{{s_i\alpha}}} 
=\left(
{1\over2}\sum_{\alpha>0} k_\alpha\langle v,s_i\alpha^\vee\rangle
t_{s_\alpha}\right)
+k_{\alpha_i}\langle v,\alpha_i^\vee\rangle t_{s_i} \cr
&=\left(
{1\over2}\sum_{\alpha>0} k_\alpha\langle s_iv,\alpha^\vee\rangle
t_{s_\alpha}\right)
+k_{\alpha_i}\langle v,\alpha_i^\vee\rangle t_{s_i} 
=
\langle s_iv,h\rangle
+k_{\alpha_i}\langle v,\alpha_i^\vee\rangle t_{s_i}. \cr
}
\formula$$
Using this equality, we obtain
$$\eqalign{
t_{s_i}[v,w]t_{s_i}
&= -t_{s_i}[\langle v,h\rangle, \langle w,h\rangle]t_{s_i} \cr
&= -[\langle s_iv,h\rangle+k_{\alpha_i}\langle v,\alpha_i^\vee\rangle
t_{s_i},
\langle s_iw,h\rangle+k_{\alpha_i}\langle w,\alpha_i^\vee\rangle
t_{s_i}] \cr
&=[s_iv,s_iw]
-k_{\alpha_i}\langle v,\alpha_i^\vee\rangle
[t_{s_i},\langle s_iw,h\rangle]
-k_{\alpha_i}\langle w,\alpha_i^\vee\rangle
[\langle s_iv,h\rangle,t_{s_i}] \cr
&=[s_iv,s_iw]
-k_{\alpha_i}\langle v,\alpha_i^\vee\rangle
(t_{s_i} \langle s_i w, h \rangle t_{s_i}
-\langle s_iw,h\rangle)t_{s_i} \cr
&\phantom{=[s_iv,s_iw]}
+k_{\alpha_i}\langle w,\alpha_i^\vee\rangle
(t_{s_i} \langle s_i v, h \rangle t_{s_i}
-\langle s_iv,h\rangle) t_{s_i} \cr
&=[s_iv,s_iw]
-k_{\alpha_i}\langle v,\alpha_i^\vee\rangle
(\langle w, h\rangle t_{s_i}
+k_{\alpha_i}\langle s_iw,\alpha_i^\vee\rangle
-\langle s_iw,h\rangle t_{s_i}) \cr
&\phantom{=[s_iv,s_iw]}
+k_{\alpha_i}\langle w,\alpha_i^\vee\rangle
(\langle v, h\rangle t_{s_i}+k_{\alpha_i}\langle s_iv,\alpha_i^\vee\rangle
-\langle s_iv,h\rangle t_{s_i}) \cr
&=[s_iv,s_iw]
-k_{\alpha_i}\langle v,\alpha_i^\vee\rangle
\langle w,\alpha_i^\vee\rangle \langle\alpha_i,h\rangle t_{s_i}
-k_{\alpha_i}^2\langle v,\alpha_i^\vee\rangle
\langle w,s_i\alpha_i^\vee\rangle  \cr
&\phantom{=[s_iv,s_iw]}
+k_{\alpha_i}\langle w,\alpha_i^\vee\rangle
\langle v,\alpha_i^\vee\rangle \langle\alpha_i,h\rangle t_{s_i}
+k_{\alpha_i}^2\langle w,\alpha_i^\vee\rangle
\langle v,s_i\alpha_i^\vee\rangle  \cr
&=[s_iv,s_iw]. \cr
}
\formula$$
The two identities (3.6) and (3.8), as in (1.3) and (1.4),
show that the algebra $A$ is isomorphic to $S(V)\otimes \CC W$.

\smallskip\noindent
(b) This can now be proved by direct computation.
If $v,w\in V$ then
$$
[\tilde v, \tilde w]
= [v-\langle v,h\rangle, w-\langle w,h\rangle] 
= [v,w] + [\langle v,h\rangle, \langle w,h\rangle]
-[v,\langle w,h\rangle]+[w,\langle v,h\rangle]=0,$$
by equation (3.4) and equation ($*$) in the proof of Theorem 3.5.
If $v\in V$ and $s_i$ is a simple reflection then, by (3.7),
$$
t_{s_i}\tilde v t_{s_i}
= t_{s_i} vt_{s_i} - 
t_{s_i} \langle v, h \rangle t_{s_i}
= s_iv - \langle s_iv,h\rangle
- k_{\alpha_i}\langle v,\alpha_i^\vee\rangle t_{s_i} 
=\widetilde{s_iv}
- k_{\alpha_i}\langle v,\alpha_i^\vee\rangle t_{s_i}. 
\qquad\hbox{\qed}
$$

Theorem 3.5b shows that if $A$ is the graded Hecke algebra defined
by (3.4), then the elements $\tilde v$, for $v\in V$, generate a subalgebra
of $A$ isomorphic to $S(V)$ and these elements together
with the $t_{s_i}$ satisfy the relations of (3.2~c).
Since part (a) of Theorem 3.5 shows that $A$ is isomorphic to $S(V)\otimes \CC W$ 
as a vector space, it follows that $A$ satisfies the conditions
(3.2~a-c), relations which uniquely define the graded Hecke algebra $A$.
Thus, Lusztig's algebras are special cases of the graded Hecke
algebras defined in Section 1.
Furthermore, by comparing the dimensions of the parameter spaces,
we see that there are graded Hecke algebras that are not isomorphic
to algebras defined by Lusztig for the Coxeter groups
$F_4$, $H_3$, $H_4$, and $I_2(m)$.

\section 4. Examples

\subsection 4A.  The symmetric group $G(1,1,n)=S_n$

Let $V$ be an $n$ dimensional vector space with orthonormal basis
$v_1,\ldots, v_n$ and let $S_n$ act on $V$ by permuting the $v_i$.
Let $A$ be a graded Hecke algebra for $S_n$.
Any element which is a product of two reflections is conjugate to
$(1,2,3)$ or $(1,2)(3,4)$.  The element $(1,2)(3,4)$ has order 2 and so,
in the algebra $A$, 
$$[v_i,v_j] = \sum_{k\ne i,j} 
(a_{(i,j,k)}(v_i,v_j)t_{(i,j,k)}+a_{(j,i,k)}(v_i,v_j)t_{(j,i,k)}),$$
since $v_i$ or $v_j$ is in $V^g=\ker a_g$ for all other three cycles $g$.
Since, by (1.6),
$a_{(j,i,k)}(v_i,v_j)=a_{(i,j,k)}(v_j,v_i)=-a_{(i,j,k)}(v_i,v_j)$,
the graded Hecke algebra $A$ is defined by the relations
$$[v_i,v_j] = \beta \sum_{k\ne i,j} (t_{(i,j,k)}-t_{(j,i,k)})
\qquad\hbox{and}\quad
t_w v_i = v_{w(i)} t_w,\formula$$
where $w\in S_n$, $1\le i,j\le n$, $i\ne j$, and 
$\beta = a_{(1,2,3)}(v_1,v_2)$.

Let $k\in \CC$.  Then, with $h$ as in (3.3),
$$
\langle v_i,h\rangle
={1\over 2}\sum_{\ell<m} k\langle v_i,v_{\ell}-v_m\rangle t_{(\ell,m)} 
={k\over 2}\left(
\sum_{i<\ell} t_{(i,\ell)} - \sum_{i>\ell} t_{(\ell,i)} \right) 
={k\over 2}\sum_{i\ne \ell} {\rm sgn}(\ell - i)t_{(i,\ell)}. 
\formula$$
If $f\in \CC S_n$, let $f\big\vert_{t_g}$ denote the coefficient of $t_g$ in $f$.
Let $A$ be the graded Hecke algebra defined by the relations in 
(4.1) with 
$$\eqalign{
\beta &= a_{(i,j,\ell)}(v_i,v_j)
=[\langle v_i,h\rangle, \langle v_j,h\rangle]\big\vert_{t_{(i,j,\ell)}} \cr
&=(k^2/4)(t_{(i,\ell)}t_{(j,\ell)}+t_{(i,j)}t_{(i,\ell)}
-t_{(j,\ell)}t_{(i,j)} )\big\vert_{t_{(i,j,\ell)}}
=k^2/4. \cr}
\formula$$
If $\tilde v_i = v_i-\langle v_i,h\rangle$ and 
$s_i$ is the simple reflection $(i,i+1)$ then, by Theorem 3.5, 
$$\eqalign{
\tilde v_i\tilde v_j&=\tilde v_j\tilde v_i,
\qquad
t_{s_i}\tilde v_i = \tilde v_{i+1}t_{s_i} + k,
\qquad
t_{s_i}\tilde v_{i+1} = \tilde v_it_{s_i} - k,
\quad\hbox{and} \cr
t_{s_j}\tilde v_i &=\tilde v_it_{s_j},
\qquad \hbox{for $|i-j|>1$,} \cr
}
\formula$$
and the algebra $A$ is the graded Hecke algebra $H_{\rm gr}$ for $S_n$ which
is defined in Section 3.  When $k=1$, the map
$$\matrix{
A &\longrightarrow &\CC S_n \cr
t_w &\longmapsto &t_w \cr
\cr
v_i &\longmapsto &{1\over 2}\displaystyle{\sum_{\ell\ne i} t_{(i,\ell)}  }\cr
}
\formula$$
is a surjective algebra homomorphism.

\subsection 4B.  The hyperoctahedral group $G(2,1,n)=WB_n$

We use the notation from Section 2B
so that the group $G(2,1,n)$ is acting by orthogonal matrices
on the $n$ dimensional vector space $V$ with orthonormal
basis $\{v_1,\ldots, v_n\}$.
In this case, $\xi_i$ denotes the diagonal matrix with 
all ones on the diagonal except for $-1$ in the $(i,i)$th entry.

Let $A$ be a graded Hecke algebra for $G(2,1,n)$.
If $\beta_1= a_{(i,j,k)}(v_i,v_j)$ and $\beta_2=a_{\xi_1(1,2)}(v_1,v_2)$,
then, in the algebra $A$,
$$
[v_i, v_j] 
=\beta_2(t_{\xi_1(1,2)}-t_{\xi_2(1,2)}) 
+\beta_1\sum_{\ell\ne i,j} 
\pmatrix{
t_{(i,j,\ell)} -t_{\xi_i\xi_\ell(i,j,\ell)}
-t_{\xi_i\xi_j(i,j,\ell)}+t_{\xi_j\xi_\ell(i,j,\ell)} \cr
+t_{\xi_i\xi_j(j,i,\ell)}+t_{\xi_j\xi_\ell(j,i,\ell)}
-t_{\xi_i\xi_\ell(j,i,\ell)}-t_{(j,i,\ell)} \cr
}.
\formula$$

Let $k_s, k_\ell\in \CC$.  Then, with $h$ as in (3.3),
$$\eqalign{
\langle v_i,h\rangle
&={k_s\over 2}\sum_\ell \langle v_i, 2v_\ell\rangle t_{\xi_\ell}
-{k_\ell\over2}\sum_{\ell<m} \langle v_i, v_\ell-v_m\rangle t_{(\ell,m)}
-{k_\ell\over2}\sum_{\ell<m} \langle v_i, v_\ell+v_m\rangle t_{\xi_\ell(\ell,m)}
\cr
&=k_s t_{\xi_i} - {k_\ell\over2}\left(
\sum_{i<\ell} (t_{(i,\ell)}+t_{\xi_i\xi_\ell(i,\ell)})
+\sum_{i>\ell} (-t_{(i,\ell)}+t_{\xi_i\xi_\ell(i,\ell)}) \right).
\cr
}
\formula$$
If $f\in \CC G(2,1,n)$, let $f\big\vert_{t_g}$ denote the coefficient of $t_g$ in $f$.
With notation as in (4.6), let $A$ be the graded Hecke algebra for $G(2,1,n)$ with
$$\eqalign{
\beta_1 &= a_{(i,j,\ell)}(v_i,v_j)
=[\langle v_i,h\rangle, \langle v_j,h\rangle]\big\vert_{t_{(i,j,\ell)}} \cr
&=(k_\ell^2/4)(t_{(i,\ell)}t_{(j,\ell)}+t_{(i,j)}t_{(i,\ell)}
-t_{(j,\ell)}t_{(i,j)} )\big\vert_{t_{(i,j,\ell)}}
=k_\ell^2/4,
\qquad\hbox{and}\cr
\beta_2 
&=[\langle v_i,h\rangle, \langle v_j,h\rangle]
\big\vert_{t_{\xi_i(i,j)}}  \cr
&=(1/2)k_s k_\ell (-t_{\xi_i}t_{(i,j)}+t_{(i,j)}t_{\xi_j}
-t_{\xi_j}t_{\xi_i\xi_j(i,j)}-t_{\xi_i\xi_j(i,j)}t_{\xi_i} )
\big\vert_{t_{\xi_i(i,j)}}
=-k_s k_\ell.
\cr
} 
$$
If $\tilde v_i = v_i-\langle v_i,h\rangle$, then, by Theorem 3.5,
the $\tilde v_i$ commute and the 
algebra $A$ is the algebra $H_{\rm gr}$ for $WB_n$ defined in Section 3.

\subsection  4C.  The type $D_n$ Weyl group $G(2,2,n)=WD_n$

We shall use the notation from Section 2B
so that the group $G(2,2,n)$ is acting by orthogonal matrices
on the $n$ dimensional vector space $V$ with orthonormal
basis $\{v_1,\ldots, v_n\}$.  This is an index 2 subgroup
of $G(2,1,n)$, and our notation is the same as used above for $WB_n$.

Let $A$ be a graded Hecke algebra for $G(2,2,n)$.
If $\beta=a_{(i,j,k)}(v_i,v_j)$ then, in the algebra $A$,
$$
[v_i, v_j] =
\beta\sum_{\ell\ne i,j} 
\pmatrix{
t_{(i,j,k)}-t_{\xi_i\xi_\ell(i,j,\ell)}
-t_{\xi_i\xi_j(i,j,\ell)}+t_{\xi_j\xi_\ell(i,j,\ell)} \cr
+t_{\xi_i\xi_j(j,i,\ell)}+t_{\xi_j\xi_\ell(j,i,\ell)}
-t_{\xi_i\xi_\ell(j,i,\ell)}-t_{(j,i,\ell)}  \cr
}.\formula$$

Let $k\in \CC$.  Then, with $h$ as in (3.3),
$$
\langle v_i,h\rangle
={k\over 2}\left(
\sum_{i<\ell} (t_{(i,\ell)}+t_{\xi_i\xi_\ell(i,\ell)})
+\sum_{i>\ell} (-t_{(i,\ell)}+t_{\xi_i\xi_\ell(i,\ell)})\right).
\formula
$$
If $f\in \CC G(2,2,n)$, let $f\big\vert_{t_g}$ denote the coefficient of $t_g$ in $f$.
With notation as in (4.8), let $A$ be the graded Hecke algebra for $G(2,2,n)$ 
with 
$$\eqalign{
\beta &= a_{(i,j,\ell)}(v_i,v_j)
=[\langle v_i,h\rangle, \langle v_j,h\rangle]\big\vert_{t_{(i,j,\ell)}} \cr
&=(k^2/4)(t_{(i,\ell)}t_{(j,\ell)}+t_{(i,j)}t_{(i,\ell)}
-t_{(j,\ell)}t_{(i,j)} )\big\vert_{t_{(i,j,\ell)}}
=k^2/4. \cr}
$$
If $\tilde v_i = v_i-\langle v_i,h\rangle$, then, by Theorem 3.5,
the $\tilde v_i$ commute and the 
algebra $A$ is the algebra $H_{\rm gr}$ for $WD_n$ defined in Section 3.

\subsection 4D.  The dihedral group $I_2(r)=G(r,r,2)$ of order $2r$

We shall use the notation for $G(r,r,2)$ from Section 2B
so that the group $G(r,r,2)$ is acting by unitary matrices
on the $2$ dimensional vector space $V$ with orthonormal
basis $\{v_1,v_2\}$.  The group $G(r,r,2)$ is realized as a 
real reflection group by using the basis
$$
\varepsilon_1 = {1\over \sqrt2}(v_1+v_2), 
\qquad
\varepsilon_2 = {-1\over i\sqrt2}(v_1-v_2).
$$
This basis is also orthonormal and, with respect to this basis, $G(r,r,2)$
acts by the matrices
$$\pmatrix{ 
\cos(2\pi m/r) &\mp\sin(2\pi m/r) \cr
\sin(2\pi m/r) &\pm\cos(2\pi m/r) \cr
}, \qquad 0\le m\le r-1.
$$

Let $A$ be a graded Hecke algebra for $G(r,r,2)$.
The conjugacy classes of elements which are products of two reflections
are $\{ \xi_1^k\xi_2^{-k}, \xi_1^{-k}\xi_2^k\}$, $0<k<r/2$.
Then, in the algebra $A$,
$$[\varepsilon_1,\varepsilon_2] = \sum_{0<k<r/2} \beta_k
(t_{\xi_1^k\xi_2^{-k}}-t_{\xi_1^{-k}\xi_2^k}),
\qquad\hbox{where}\quad
\beta_k=a_{\xi_1^k\xi_2^{-k}}(\varepsilon_1,\varepsilon_2).
\formula$$

When $r$ is even, there are two conjugacy classes of reflections
$$\{\xi_1^{2k}\xi_2^{-2k}(1,2) \ |\ 0\le k<r/2\}
\qquad\hbox{and}\qquad
\{\xi_1^{2k+1}\xi_2^{-(2k+1)}(1,2) \ |\ 0\le k<r/2\}.$$
The reflection $\xi_1^m\xi_2^{-m}(12)$ is the reflection in the 
line perpendicular to the vector 
$$\alpha_m = \sin(-2\pi m/2r )\varepsilon_1+\cos(-2\pi m/2r )\varepsilon_2,$$
and the vectors $\alpha_m$ can be taken as a root system
for $G(r,r,2)$.  With $h$ as in (3.3) and $k_s, k_\ell \in \CC$, 
$$\eqalign{
\langle \varepsilon_1,h\rangle 
&= \sum_{0\le k<r/2} \left(k_s\sin(-2k\, 2\pi/2r)t_{\xi_1^{2k}\xi_2^{-2k}(1,2)}
+k_\ell\sin(-(2k+1)\, 2\pi/2r)t_{\xi_1^{2k+1}\xi_2^{-(2k+1)}(1,2)}\right), \cr
\langle \varepsilon_2,h\rangle 
&= \sum_{0\le k<r/2} \left(k_s\cos(-2k\, 2\pi/2r)t_{\xi_1^{2k}\xi_2^{-2k}(1,2)}
+k_\ell\cos(-(2k+1)\, 2\pi/2r)t_{\xi_1^{2k+1}\xi_2^{-(2k+1)}(1,2)}\right). \cr
}
\formula$$
If $f\in \CC G(r,r,2)$, let $f\big\vert_{t_g}$ denote the coefficient of $t_g$ in $f$.
With notation as in 4.10, let $A$ be the graded Hecke algebra for
$G(r,r,2)$ with
$$
\beta_k = a_{\xi_1^k\xi_2^{-k}}(\varepsilon_1,\varepsilon_2)
=[\langle \varepsilon_1,h\rangle, 
\langle \varepsilon_2,h\rangle]\big\vert_{t_{\xi_1^k\xi_2^{-k}}}
=\cases{
\sin(k \, 2\pi /2r) r k_sk_\ell &if $k$ is odd \cr
\cr
\sin(k \, 2\pi  /2r) {r\over 2}(k_s^2+k_\ell^2) &if $k$ is even. \cr
}
\formula
$$
If $\tilde \varepsilon_i = \varepsilon_i-\langle \varepsilon_i,h\rangle$, 
then by Theorem 3.5,
the $\tilde \varepsilon_i$ commute and the 
algebra $A$ is the algebra $H_{\rm gr}$ for $I_2(r)$ defined in Section 3.

When $r$ is odd, all aspects of the calculation 
in (4.11) and (4.12) are the same as for the case $r$ even
except that there is only one conjugacy class of reflections,
$\{\xi_1^k\xi_2^{-k}(1,2) \ |\ 0\le k\le r-1\}$,
and so $k_s=k_\ell$.

\subsection 4E. The group $G(r,r/2,2)$, $r/2$ odd

We use the notation from Section 2B, or from above for the group $G(r,r,2)$.
In this case, the group is not a real reflection group, hence 
$G(r,r/2,2)$ acts by unitary matrices but not by orthogonal matrices.

Let $A$ be a graded Hecke algebra for $G(r,r/2,2)$.
The only conjugacy class for which $a_g$ can be nonzero is 
$\{t_{\xi_1^k\xi_2^{r/2-k}(1,2)}\ |\ 0\le k<r\}$.  Thus, in the algebra $A$,
$$[v_1,v_2] = \beta\sum_k (t_{\xi_1^{2k}\xi_2^{r/2-2k}(1,2)}
-t_{\xi_1^{r/2-2k}\xi_2^{2k}(1,2)}),
\qquad\hbox{where}\quad
\beta = a_{\xi_2^{r/2}(1,2)}(v_1,v_2).$$

\vfill\eject

\section 5.  A different graded Hecke algebra for $G(r,1,n)$

The classification of graded Hecke algebras for complex 
reflection groups in Section 2 shows that
there do not exist graded Hecke algebras $A\cong S(V)\otimes \CC G$
for the groups $G=G(r,1,n)$, $r>2$, $n>3$.
In this section, we define a different ``semidirect product'' of
the symmetric algebra $S(V)$ and the group algebra $\CC G$ 
for the groups $G(r,1,n)$.  These algebras are not graded Hecke
algebras in the sense of Section 1, but they do have a structure similar
to what we would expect from experience with graded Hecke algebras
for real reflection groups.
Is it possible that there is a general definition of graded Hecke
algebras, different from that given in Section 1, which includes
the algebras defined below as examples for the groups $G(r,1,n)$?

We shall use the notation for the groups $G(r,1,n)$ as in Section 2B
so that the group $G(r,1,n)$ is acting by monomial matrices on
a vector space $V$ of dimension $n$ with orthonormal basis
$\{v_1,\ldots, v_n\}$.  Let $s_i$ denote the permutation $(i,i+1)\in G(r,1,n)$.

Define $H^*_{r,1,n}$ to be the algebra generated by
the group algebra $\CC G(r,1,n)$ and $V$
with relations
$$\matrix{
\displaystyle{
v_iv_j=v_jv_i },\hfill &\qquad &\hbox{for all $1\le i,j\le n$,} \hfill \cr
\cr
\displaystyle{
t_{\xi_i}v_j=v_jt_{\xi_i} },\hfill &\qquad &\hbox{for all $1\le i,j\le n$,} \hfill \cr
\cr
\displaystyle{
t_{s_i}v_k=v_kt_{s_i} },\hfill &\qquad &\hbox{if $k \notin \{i, i+1 \}$,} \hfill\cr
\cr
\displaystyle{
t_{s_i}v_{i+1}=v_it_{s_i}+\sum_{\ell=0}^{r-1} t_{\xi_i^\ell\xi_{i+1}^{-\ell}} }\,,
\hfill
&&\hbox{for $1\le i\le n-1$.}\hfill \cr
}
\formula
$$
The following proposition establishes an ``evaluation homomorphism''
for the algebras $H_{r,1,n}^*$ which is a generalization of the homomorphism
in (4.5).

\prop  Define elements $\bar v_k$ in the group algebra $\CC G(r,1,n)$ by
setting $\bar v_1=0$ and 
$$\bar v_k = {1\over r}\ \sum_{i<k} \ \sum_{0\le \ell\le r-1} 
t_{\xi_i^\ell\xi_k^{-\ell}(i,k)},
\qquad\hbox{for $2\le k\le n$.}
$$
Then there is a surjective algebra homomorphism
$$\matrix{
H^*_{r,1,n} &\longrightarrow & \CC G(r,1,n) \cr
t_g &\longmapsto &t_g \cr
v_k &\longmapsto &\bar v_k \cr}
$$
\pf
We must check that the defining relations (5.1) of $H^*_{r,1,n}$ hold 
with the $v_k$ replaced by the $\bar v_k$.

For each $1\le k\le n$, let
$$\bar z_k = \bar v_1+\cdots+\bar v_k
= {1\over r}\sum_{1\le i<j\le k\atop 0\le \ell\le r-1} t_{\xi_i^\ell\xi_j^{-\ell}(i,j)}.$$
Then, for each $k$, $\bar z_k\in Z(\CC G(r,1,k))$ since it is the sum of the
elements of the conjugacy class of reflections $t_{\xi_i^\ell\xi_j^{-\ell}(i,j)}$
in $G(r,1,k)$.  So $\bar z_k$ commutes with $\bar z_1,\ldots,\bar z_k$ and
therefore $\bar z_1,\ldots, \bar z_n$ commute.  Since
$\bar v_k=\bar z_k-\bar z_{k-1}$, it follows that $\bar v_1,\ldots, \bar v_n$
also commute.

If $m>k$ then $t_{\xi_m}$ clearly commutes with $\bar z_k$.  If $m\le k$ then
$t_{\xi_m}$ commutes with $\bar z_k$ since $\bar z_k\in Z(G(r,1,k))$.  So
$t_{\xi_m}$ commutes with $\bar z_1,\ldots,\bar z_n$ and hence
with $\bar v_1,\ldots,\bar v_n$.

Since
$$\eqalign{
t_{s_k}\bar v_kt_{s_k}
\ & = \ t_{s_k}\left(\sum_{i<k\atop 0\le \ell\le r-1} 
t_{\xi_i^\ell\xi_k^{-\ell}(i,k) }\right) t_{s_k} 
\ = \ \sum_{i<k\atop 0\le \ell\le r-1} 
t_{\xi_i^\ell\xi_{k+1}^{-\ell}(i,k+1)} \cr
&= 
{\vrule height 25pt depth6pt width0pt}
\sum_{i<k+1\atop 0\le \ell\le r-1} 
t_{\xi_i^\ell\xi_{k+1}^{-\ell}(i,k+1) }
\ -\ \sum_{0\le \ell\le r-1} t_{\xi_k^\ell\xi_{k+1}^{-\ell}(k,k+1)} \cr
&= 
{\vrule height 25pt depth6pt width0pt}
\ \ \ \bar v_{k+1} \ -
\sum_{0\le \ell\le r-1} t_{\xi_k^\ell\xi_{k+1}^{-\ell}}\,
t_{s_k}, \cr
}
$$
it follows that
%$$t_{s_k}\bar v_k = \bar v_{k+1} t_{s_k} 
%-\sum_{0\le \ell\le r-1} t_{\xi_k^\ell\xi_{k+1}^{-\ell}}.$$
%Multiplying left and right by $t_{s_k}$ gives
$$
\bar v_k t_{s_k}
=\ t_{s_k} \bar v_{k+1}  \
-\sum_{0\le \ell\le r-1} t_{s_k} t_{\xi_k^{\ell}\xi_{k+1}^{-\ell}} t_{s_k} 
=\ t_{s_k} \bar v_{k+1} \ 
-\sum_{0\le \ell\le r-1} t_{\xi_k^{-\ell}\xi_{k+1}^\ell} 
=\ t_{s_k} \bar v_{k+1} \ 
-\sum_{0\le \ell\le r-1} t_{\xi_k^{\ell}\xi_{k+1}^{-\ell}}. 
$$
\hfill
\qquad\hbox{\qed}

%\vfill\eject

\section 6. References

\medskip\noindent
\item{[Bou]} {\smallcaps N.\ Bourbaki},
{\sl Groupes et algebres de Lie}, 
Chapt.\ IV-VI, Masson, Paris, 1981.

\medskip\noindent
\item{[Ca]} {\smallcaps R.W.\ Carter},
{\it Conjugacy classes in the Weyl group}, Compositio Math.\ {\bf 25} (1972), 1-59.

\medskip\noindent
\item{[Dr]} {\smallcaps V.G.\ Drinfeld},
{\it Degenerate affine Hecke algebras and Yangians},
Funct.\ Anal.\ Appl.\ {\bf 20} (1986), 58-60.

\medskip\noindent
\item{[G+]}
{\smallcaps
M.\ Geck, G.\ Hiss, F.\ L\"ubeck, G.\ Malle, and G.\ Pfeiffer},\
{\it CHEVIE -- A system for computing and processing 
generic character tables for finite groups
of Lie type, Weyl groups and Hecke algebras}, 
AAECC, {\bf 7}  (1996), 175-210. 

\medskip\noindent
\item{[Hu]}
{\smallcaps J.E.\  Humphreys}, 
{\sl Reflection groups and Coxeter groups}, Cambridge Studies in Advanced Mathematics
{\bf 29}, Cambridge University Press, Cambridge, 1990.

\medskip\noindent
\item{[Lu]} {\smallcaps G.\ Lusztig},
{\it Affine Hecke algebras and their graded version},
J.\ Amer.\ Math.\ Soc.\ {\bf 2} (1989), 599-635.

\medskip\noindent
\item{[Lu2]} {\smallcaps G.\ Lusztig},
{\it Cuspidal local systems and graded Hecke algebras I},
Inst.\ Hautes \'Etudes Sci.\ Publ.\ Math.\ {\bf 67} (1988), 145-202.

\medskip\noindent
\item{[Mac]} {\smallcaps I.G.\ Macdonald},
{\sl Symmetric functions and Hall polynomials},
Second edition, Oxford Mathematical Monographs, Oxford University Press, 
New York, 1995.

\medskip\noindent
\item{[OS1]} {\smallcaps P.\ Orlik and L.\ Solomon},
{\it Coxeter arrangements}, 
Proc.\ Symp.\ Pure Math.\ {\bf 40} (1983), 269-291.

\medskip\noindent
\item{[OS2]} {\smallcaps P.\ Orlik and L.\ Solomon},
{\it Arrangements defined by unitary reflection groups},
Math.\ Ann.\ {\bf 261} (1982), 339-357.

\medskip\noindent
\item{[OT]} {\smallcaps P.\ Orlik and H.\ Terao},
{\it Arrangements of hyperplanes},
Springer-Verlag, Berlin-Heidelberg, 1992.

\medskip\noindent
\item{[S+]}
{\smallcaps Martin Schönert et.al.},
{\it GAP -- Groups, Algorithms, and Programming},
Lehrstuhl D f\"ur Mathematik, Rheinisch Westf\"alische Technische
Hochschule, Aachen, Germany, fifth edition, 1995.

\medskip\noindent
\item{[ST]} {\smallcaps G.C.\ Shephard and J.A.\ Todd},
{\it Finite Unitary Reflection Groups}, Canad.\ J.\ Math.\ 
{\bf 6} (1954), 274-304.

\medskip\noindent
\item{[St1]} {\smallcaps R.\ Steinberg},
{\it Differential equations invariant under finite reflection
groups}, Trans.\ Amer.\ Math.\ Soc.\ {\bf 112} (1964), 392-400.

\medskip\noindent
\item{[St2]} {\smallcaps R.\ Steinberg},
{\it Endomorphisms of linear algebraic groups},
Mem.\ Amer.\ Math.\ Soc.\ {\bf 80} (1968), 1-108.

\vfill\eject

\end